\documentclass{article}
\usepackage{amsfonts}

\usepackage{graphics}

\usepackage{amsmath}

\usepackage{amssymb}

\usepackage{amscd}

\usepackage{euscript}

\usepackage{color}

\newtheorem{Theorem}{Theorem}

\newtheorem{Corollary}[Theorem]{Corollary}

\newtheorem{Example}[Theorem]{Example}

\newtheorem{Definition}[Theorem]{Definition}

\newtheorem{Lemma}[Theorem]{Lemma}

\newtheorem{Proposition}[Theorem]{Proposition}

\newtheorem{Fundamental Theorem}{Fundamental Theorem}

\newenvironment{Proof}[1][Proof]{\textbf{#1.} }{\ \rule{0.5em}{0.5em}}

\def \la {\mathfrak{a}}
\def \lb {\mathfrak{b}}
\def \ad {\mathrm{ad}}
\def \A {\EuScript{A}}
\def \Ad {\mathrm{Ad}}
\def \a {\alpha}
\def \Aut {{\rm Aut}}

\def \C {\mathcal{C}}

\def \Der {\mathrm{Der}}

\def \ga {\gamma}

\def \e {\epsilon}

\def \w {\omega}

\def \id {\mathrm{id}}

\def \lg {\mathfrak {g}}

\def \R {\mathbb{R}}

\def \X {\EuScript{X}}

\def \W {\Omega}

\def \G {\Gamma}

\def \f {\phi}

\def \d {\partial}

\def \M {\EuScript{M}}

\def \N {\mathbb N}

\def \n {\mu}

\def \Rank {\mathrm{Rank}}

\def \P {\EuScript{P}}

\def \a {\alpha}

\def \b {\beta}

\def \g {\gamma}

\def \W {\Omega}

\def \t {\widetilde}

\def \tn {\otimes}

\def \tr {\triangleright}

\def \Gc {\mathcal{G}}
\def \L {{\mathcal{L}}}
\def \ra {\xrightarrow}

\def \le {\mathfrak{e}}

\def \S {{\EuScript{S}}}

\def \Q {{\EuScript{Q}}}

\def \Qc {{\mathcal{Q}}}

\def \Wc {{\EuScript{W}}}

\def \diff {\mathrm{diff}}

\def \Sc {\mathcal{S}}

\def \Pc {\mathcal{P}}

\def \H {\mathcal{H}}

\def \Z {\mathbb{Z}}

\begin{document}

\author{Jo\~{a}o Faria Martins  \\{\small Centro de Matem\'{a}tica da Universidade do Porto}\\{\small
Rua do Campo Alegre, 687, 4169-007 Porto, Portugal} \\ {\it \small jmartins@math.ist.utl.pt}\\ \\ Roger Picken  \\ {\small Departamento de Matem\'{a}tica,} \\ {\small Instituto Superior T\'{e}cnico, TU Lisbon } \\ {\small Av. Rovisco Pais, 1049-001 Lisboa, Portugal} \\ {\it \small rpicken@math.ist.utl.pt}}

\title{On Two-Dimensional Holonomy}

\maketitle

\begin{abstract}
We define the thin fundamental categorical group  ${\mathcal P}_2(M,*)$ of a based smooth manifold  $(M,*)$ as the categorical group whose objects are rank-1 homotopy classes of based loops on $M$, and whose morphisms are rank-2 homotopy classes of homotopies between based loops on $M$. Here two maps are rank-$n$ homotopic, when the rank of the differential of the homotopy between them equals $n$. Let $\C(\Gc)$ be a Lie categorical group coming from a Lie crossed module ${\Gc= ( \d\colon E \to  G,\tr)}$. We construct categorical holonomies, defined to be smooth morphisms ${\mathcal P}_2(M,*) \to \C(\Gc)$,  by using a notion of categorical connections, being a pair $(\w,m)$, where $\w$ is a connection 1-form on $P$, a principal $G$ bundle over $M$, 
and $m$ is a 2-form on $P$ with values in the Lie algebra of $E$, with the pair $(\w,m)$ satisfying suitable conditions.
 As a further result, we are able to define Wilson spheres in this context.
\end{abstract}


\noindent{\bf 2000 Mathematics Subject Classification} {\it 53C29 
(primary); 
18D05 (secondary) 
}

\newpage
\tableofcontents

\section*{Introduction}

Categorification is an influential idea in many areas of mathematics, and in geometry it is natural to try and think about categorifying the notions of holonomy and parallel transport in terms of higher categorical generalisations of the notions of loop,  Lie group and connection on a principal bundle, in the spirit of Baez and Schreiber \cite{B, BS}.
 In this article we construct a framework for $2$-dimensional, or surface, holonomy along such lines. The based loops on a manifold $M$ are replaced by what we call the (strict) thin fundamental categorical group of $M$, ${\mathcal P}_{{2}}(M,*)$, a monoidal category whose objects are rank-1 homotopy classes of based loops on $M$ and whose morphisms are rank-2 homotopy classes of homotopies between based loops (or $1$-parameter families of loops).  Here, rank-1 homotopy, at the level of loops, means that the loops are homotopic in such a way that the 
rank of the differential of the homotopy between them is less than or equal to $1$, i.e. 
the homotopies between loops do not sweep out area. Similarly, at the level of morphisms, two homotopies between based loops are rank-2 homotopic, if they themselves are homotopic in such a way that the 
rank of the differential of the homotopy between them is less than or equal to $2$, i.e. 
the homotopy between homotopies does not sweep out volume. For precise definitions, we refer to subsections \ref{1-Tracks} and \ref{2-Tracks}. The Lie group is replaced by a categorical Lie group, $\C(\Gc)$, naturally obtained from a crossed module of Lie groups ${\Gc= ( \d\colon E \to  G,\tr)}$. The connection on a principal $G$-bundle $P$ over $M$ is replaced by what we call a categorical connection, consisting of a {connection} $1$-form on $P$ with values in $\lg$, the Lie algebra of $G$, together with a $2$-form on $P$ with values in $\le$, the Lie algebra of $E$, satisfying some conditions including the well-known ``vanishing of the fake curvature''. Locally these conditions correspond to Baez and Schreiber's \cite{BS} local formulation for a $2$-connection. The Maurer-Cartan structure equation and Bianchi identity for the curvature of an ordinary connection are shown to have natural analogues for the curvature $3$-form of the categorical connection.  Also a variant of the Ambrose-Singer theorem (which plays a crucial role in our construction of categorical holonomy) generalises to a higher-order version.  

A categorical holonomy is defined to be a (strict monoidal) {smooth} functor from ${\mathcal P}_{{2}}(M,*)$ to  $\C(\Gc)$. The main result that we prove (Thm. \ref{FFF}) is how a categorical connection gives rise to a categorical holonomy. The underlying geometrical idea is to lift the $1$-parameter family of loops into $P$,  horizontally in one direction, namely the direction along the loops, and to use this lift to pull back the forms of the categorical connection, which are then integrated suitably.

Note that the appearance of a principal $G$-bundle $P$ with a connection is natural in the context of categorical holonomies. This is because, at the level of the set of objects of ${\mathcal P}_{{2}}(M,*)$, i.e.  $\pi_1^1(M,*)$, the thin homotopy classes of based loops on $M$,  any smooth functor  ${\mathcal P}_{{2}}(M,*) \to \Gc$ gives rise to a smooth group morphism $\pi_1^1(M,*) \to G$, and therefore it defines a principal $G$-bundle over $M$ with a connection \cite{CP,MP}.

The whole construction is carried out using the language and methods of differential geometry and principal bundles, thereby avoiding working with infinite-dimensional path spaces, which was {an} approach taken in \cite{BS}.  The construction is coordinate-free from the outset, since we use forms defined on $P$. We remark that the construction of ${\mathcal P}_{{2}}(M,*)$ is of interest in its own right in defining a strict thin fundamental categorical group  of a manifold (previously only a weak version was known). This part of our construction is very similar to that of \cite{BHKP,HKK}.

 {We study the relation between our construction and non-abelian gerbes, {as in \cite{BrMe},} and 2-bundles, as in \cite{BS}, in  \ref{2bundles}. Note that each 2-bundle with structure 2-group coming from a crossed module of the form  $(\Ad\colon G \to \Aut(G),\tr)$, where $\tr$ denotes the obvious left action of $\Aut(G)$ in $G$, is naturally a non-abelian gerbe.}  Let $\Gc=(\d\colon E \to G,\tr)$ be a Lie crossed module. Our construction corresponds to a particular case of $\Gc$-2-bundles, for which the $E$-valued transition functions are trivial, and therefore the $G$-valued transition functions  satisfy the usual cocycle condition for a principal $G$-bundle.
 Although $\Gc$-2-bundles are a natural way to approach two-dimensional holonomy, we emphasise that our main goal is the definition of categorical group maps ${\mathcal P}_{{2}}(M,*) \to \C(\Gc)$. {This is a very strong condition, and
we argue that it should force the {${\mathcal G}$-2-bundle} with connection to be of
the special form considered here - see 2.4.6 where this point is elaborated.}

{At the end of this article we define the notion of Wilson sphere, which  means that a categorical holonomy, taking values in the kernel of $\d\colon E\to G$, can be associated to embedded spheres in a manifold, up to acting by an element of the group $G$. }

\section{Preliminaries}

\subsection{Resum\'{e} of results and notation from the theory of principal fibre bundles}

\subsubsection{Connections on principal fibre bundles}\label{CCC}

For full details on connections on principal fibre bundles see \cite{KN,P}.

Let $M$ be a paracompact smooth manifold. Let $G$ be a Lie group with Lie algebra $\lg$.  Consider a principal $G$-bundle $P$ over $M$, with canonical projection $\pi\colon P \to M$. For each $x \in M$, let $P_x=\pi^{-1}(x)$ be the fibre at $x$. Let also $\X(P)$ denote the Lie algebra  of smooth vector fields on $P$. 
Recall that there exists a Lie algebra morphism $A \in \lg \mapsto A^\# \in \X(P)$, such that $A^\#_u=\frac{d}{dt}u \exp(tA)_{t=0}, \forall u \in P, \forall A \in \lg$.  Here if $X\in \X(P)$ is a vector field then $X_u \in T_u P$ (the tangent space of $P$ at $u\in P$) denotes the value of $X$ at the point $u$.

Denote the right action of $G$ on $P$ as $g \in G \mapsto R_g \in \mathrm{diff}(P)$, {where ${\rm diff}(P)$ denotes the diffeomorphism group of $P$}.  Recall the following very useful formula: 
$${(R_g)}_*(A^\#)=\left (g^{-1}Ag\right )^\#,$$ where $g \in G$ and $A \in \lg$.
Given an element $u \in P$, recall that a vector $X_u \in T_u P$  is said to be vertical if $\pi_*(X_u)  =0$. Let $T_u^VP$ denote the subspace of vertical vectors in  $T_u P$. Notice that any $X_u \in T_u^VP$ is of the form $A^\#_u$ for some $A \in \lg$, and this correspondence is one-to-one.

 Let $\w\in \A^1(P,\lg)$ be a connection  on $P$. In other words $\w$ is a smooth 1-form on $P$ with values in $\lg$ such that:
\begin{enumerate}
\item $R_g^*(\w)=g^{-1}\w g , \forall g \in G,$
\item $\w(A^\#)=A,\forall A \in \lg$.
\end{enumerate}
Given an element $u \in P$, let $T^H_u P=\{X_u \in T_u P\colon \w(X_u)=0\}$, the horizontal subspace at $u \in P$. Then $T^H_u$ is a complementary subspace of $T^V_uP$ in $T_u P$ for each $u \in P$. Moreover, the following identity holds for each $g \in G$ and $u \in P$:
$$\left (T^H_u P\right ) g= T^H_{ug} P .$$

The natural projection maps $X \in TP \mapsto X^H$ and $X \in TP \mapsto X^V$ are smooth.  Here $TP$ denotes the tangent bundle of $P$.
There exists also a unique map     
$X \in \X(M) \mapsto \t{X}\in \X(P)$, called the horizontal lift of $X\in \X(M)$, such that $\t{X}_u \in T_u^HP, \forall u \in P$ and $\pi_*(\t{X})=X$, for any $X \in \X(M)$, the Lie algebra of smooth vector fields on $M$.  Note that this horizontal lift {gives rise to linear maps} $T_x M  \to T_u P$, if $\pi(u)=x$, where $x \in M$ and $u \in P$.  The horizontal lift of a vector field $X \in \X(M)$ is always $G$-{invariant}. In other words:
$$\t{X} g=\t{X}, \forall g \in G.$$

\subsubsection{Curvature}\label{curv}

Let $P$ be a principal $G$-bundle over $M$. Let $\w \in \A^1(P,\lg)$ be a connection 1-form on $P$. Given an $n$-form $a$ 
on $P$, the exterior covariant derivative of $a$ is given by $$Da=d a \circ (\underbrace{H \times H \ldots \times  H}_{(n+1)\textrm{-times }} ).$$

Let $\W\in \A^2(P,\lg)$ be the curvature 2-form of the connection $\w$. It can be defined as the exterior covariant derivative $D\w$ of the connection 1-form $\w$, in other words:
$$\W(X,Y)=d \w(X^H,Y^H),$$
where $X,Y \in \X(P)$. 
{The curvature $2$-form $\W$ is $G$-{equivariant}, which means:} $$R_g^*(\W)=g^{-1}\W g, \forall g \in G .$$ 
  Recall also Cartan's structure equation:
\begin{equation}\label{cse}{d \w(X,Y) +[\w(X),\w(Y)]=\W(X,Y),}\end{equation}
valid for any vector fields $X,Y \in \X(P)$. The Bianchi identity can be stated by saying that the exterior covariant derivative $D \W$ of the curvature form is zero, in other words:
\begin{equation}D \W=d \W\circ (H \times H \times H)=0.\end{equation}
{The Bianchi identity can also be written as:}
$${d \W(X,Y,Z)+[\w(X),\W(Y,Z)]+[\w(Y),\W(Z,X)]+[\w(Z),\W(X,Y)]=0,}$$
{for any smooth vector fields $X,Y,Z \in \X(P)$.}

\subsubsection{{Parallel} transport}\label{Hol}

Let $P$ be a principal $G$-bundle over the manifold $M$. Let $\w \in \A^1(P,\lg)$ be a connection on $P$.  Recall that $\w$ determines a parallel {transport} along smooth curves. Specifically, given $x \in M$ and a smooth curve $\gamma\colon [0,1] \to M$, with $\gamma(0)=x$, then there exists a smooth map:
 $$ (t,u) \in [0,1] \times P_x \mapsto {\mathcal{H}}_\w(\ga,t,u) \in  P,$$ 
uniquely defined by  the conditions:

\begin{enumerate} 
\item $\frac{d}{d t} {\mathcal{H}}_\w(\ga,t,u)=\left(\t{\frac{d}{dt} \ga(t)}\right)_{{\mathcal{H}}_\w(\ga,t,u) };\forall t \in [0,1], \forall u \in P_x,$
\item ${\mathcal{H}}_\w(\ga,0,u)=u; \forall u \in P_x. $
\end{enumerate}
In particular this implies that ${\mathcal{H}}_\w(\ga,t)$, {given by $u\mapsto \mathcal{H}_\w(\ga,t,u)$}, maps $P_x$ bijectively into  $P_{\ga(t)}$, for any $t \in [0,1]$. 
Recall that the parallel transport is $G$-{equivariant}, in other words:
$$\H_\w(\ga,t,ug)=\H_\w(\ga,t,u)g, \forall g \in G, \forall u \in P_x.$$

\subsubsection{The dependence of the parallel transport on a smooth family of curves - the Ambrose-Singer theorem}\label{BBB}

Let $M$ be a smooth manifold. Let ${{\rm D}}^n\doteq [0,1]^n$ be the $n$-cube, where $n \in \N$. A map $f\colon  {{\rm D}}^n \to M$ is said to be smooth if {its partial derivatives of any order} exist and are continuous as maps ${{\rm D}}^n \to M$. 
Let $G$ be a Lie group with Lie algebra $\lg$. Consider a smooth principal $G$-bundle $\pi \colon P \to M$ with  a connection 1-form $\w \in \A^1(P,\lg)$.

Now let  $s \in [0,1] \mapsto \g_s$ be a smooth 1-parameter family of smooth curves $[0,1] \to M$. Here smooth means that the map $\G\colon (t,s) \in [0,1]^2 \mapsto \ga_s(t) \in M$ is smooth. Define the initial point map $q\colon [0,1] \to M$ by $q(s)=\g_s(0)$ for each $s\in[0,1]$. Choose $u \in P_{q(0)}$ and set $u_s={\mathcal{H}}_\w(q,s,u)\in P_{q(s)}$, where $s \in  [0,1]$. Our purpose is to analyse the $s$-dependence of  ${\mathcal{H}}_\w(\ga_s,t,u_s)$, where $s,t \in  [0,1]$, by calculating:
\begin{align*}
\frac{\d}{\d s} {\mathcal{H}}_\w(\g_s,t,u_s).
\end{align*}
This analysis is of course classical.

For convenience of notation, define
 $$ f(s,t)={\mathcal{H}}_\w(\g_s,t,u_s);s,t \in [0,1].$$
 Let the coordinate vector fields of $[0,1]^2$ be $\frac{\d}{\d t}$ and $\frac{\d}{\d s}$.
We have:
\begin{align*}
\frac{\d}{\d t} \w \left (\frac{\d}{\d s} f(s,t) \right) & =  \frac{\d}{\d t} f^*(\w) \left (  \frac{\d}{\d s} \right )\\
&=f^*(d \w) \left ( \frac{\d}{\d t}, \frac {\d}{\d s}\right)+\frac{\d}{\d s} f^*(\w) \left (  \frac{\d}{\d t} \right )\\
&=d \w \left ( \frac{\d}{\d t}f(s,t), \frac {\d}{\d s}f(s,t)\right).
\end{align*}
Note that the second {equation follows from}  the well-known formula
\begin{equation}\label{A}
d\f(X,Y)=X\f(Y)-Y\f(X)-\f([X,Y]), 
\end{equation}
valid for any 2-form $\f$  on a manifold $M$ and any $X,Y \in \X(M)$.
The last {equation} follows from the fact that $\frac{\d}{\d t} f(s,t)$ is horizontal, by definition of parallel transport.  

Therefore by Cartan's structure equation (\ref{cse}) and the fact that $\frac{\d}{\d t} f(s,t)$ is horizontal it follows that:
\begin{equation}\label{CC}
\frac{\d}{\d t} \w \left (\frac{\d}{\d s} f(s,t) \right) =\W\left ( \frac{\d}{\d t}f(s,t), \frac {\d}{\d s}f(s,t)\right).
\end{equation}
Since  $\W(X,Y)=0$ if either of the vectors $X$ or $Y$  is vertical we can conclude:
$$
\frac{\d}{\d t} \w \left (\frac{\d}{\d s} f(s,t) \right) =\W\left ( \left (\frac{\d}{\d t}f(s,t)\right )^H, \left (\frac {\d}{\d s}f(s,t)\right)^H\right).$$

 Note that, by definition of parallel transport, we have  $$\left (\frac{\d}{\d t}f(s,t)\right)^H=\frac{\d}{\d t}f(s,t)= \left (\t{\frac{\d}{\d t}\g_s(t)}\right)_{{\mathcal{H}}_\w(\g_s,t,u_s) }.$$
Given that $\left (\frac{\d}{\d s}f(s,t)\right )^H$ is horizontal and  the fact that:
$$
\pi_*\left (\left (\frac{\d}{\d s}f(s,t)\right)^H \right) =\pi_* \left(\frac{\d}{\d s}f(s,t) \right)=\frac{\d}{\d s}\pi \big ( f(s,t)\big )=\frac{\d}{\d s}\g_s(t),
$$
it follows
$$\left (\frac{\d}{\d s}f(s,t)\right)^H= \left (\t{\frac{\d}{\d s}\g_s(t)}\right)_{{\mathcal{H}}_\w(\g_s,t,u_s) }.$$

{Going back to equation (\ref{CC}), we have that, {for each $s \in [0,1]$,} }
$\w \left (\frac{\d}{\d s} f(s,0) \right)=0,$
since $f(s,0)=u_s={\mathcal{H}}_\w(q,s,u)$, hence $\frac{\d}{\d s} f(s,0)$ is horizontal for each $s \in [0,1]$.  We thus arrive at the following {well known} {result}:

\begin{Lemma}\label{Main1}
Let $G$ be a Lie group with Lie algebra $\lg$. Let $P$ be a smooth principal $G$-bundle over the manifold $M$. Let $s \in [0,1] \mapsto \g_s$ be a smooth 1-parameter family of curves $[0,1] \to M$. Here smooth means that the map $(s,t) \in [0,1]^2 \mapsto \g_s(t)\in M$ is smooth. Consider a connection $\w \in \A^1(P,\lg)$.  Let $q\colon [0,1] \to M$ be the curve such that $q(s)=\g_s(0), \forall s \in [0,1]$.
Choose $u \in P_{q(0)}$, and let  $u_s={\mathcal{H}}_\w(q,s,u)\in P_{q(s)}$, where $s \in  [0,1]$.  The following holds for each $s,t' \in [0,1]$:
\begin{equation}
\w\left (\frac{\d}{\d s} { {\mathcal{H}}_\w(\g_s,t',u_s)}\right)=\int_{0}^{t'} \W \left (\t{\frac{\d}{\d t}\g_s(t)},\t{\frac{\d}{\d s}\g_s(t)} \right)_{{\mathcal{H}}_\w(\g_s,t,u_s) }  d t.
\end{equation}
\end{Lemma}

Consider  the case when all curves  $\g_s\colon [0,1] \to M$, where  $s \in [0,1]$, are closed, though possibly with varying initial points. (This flexibility will be important to define Wilson spheres in \ref{WS}). 
Then we can define a family of holonomies $g_{\g_s}\in G$, where $s \in [0,1]$, by:
 $$u_sg_{\g_s}= {\mathcal{H}}_\w(\g_s,1,u_s), \forall s \in [0,1].$$
We will use the preceding lemma to obtain a differential equation satisfied by $g_{\g_s}$.

Because $\frac{d}{d s} u_s$ is horizontal it follows that:
$$
{\w \left ( \frac{d}{d s} (u_s g_{\g_s})\right) = \w \left ( u_s\frac{d}{d s} g_{\g_s}\right)=\int_{0}^{1} \W \left (\t{\frac{\d}{\d t}\g_s(t)},\t{\frac{\d}{\d s}\g_s(t)} \right)_{{\mathcal{H}}_\w(\g_s,t,u_s) }dt.}
$$
{Since, for each $s \in [0,1]$, the vector $u_s\frac{d}{d s} g_{\g_s}$ is vertical, this means,   by the second condition of the definition of a connection 1-form $\w$ (see \ref{CCC}):}
\begin{equation}\label{D}
u_s\frac{d}{d s} g_{\g_s} =\left (\int_{0}^{1} \W \left (\t{\frac{\d}{\d t}\g_s(t)},\t{\frac{\d}{\d s}\g_s(t)} \right)_{{\mathcal{H}}_\w(\g_s,t,u_s)  } dt\right)^\#_{{\mathcal{H}}_\w(\g_s,1,u_s) },
\end{equation}
or
\begin{align*}
&\left [\frac{d}{d t'}u_s g_{\g_{(s+t')}}\right]_{t'=0}\\&\quad\quad \quad \quad=\frac{d}{d t'} {{\mathcal{H}}_\w(\g_s,1,u_s) }\exp\left (t'\int_{0}^{1} \W \left (\t{\frac{\d}{\d t}\g_s(t)},\t{\frac{\d}{\d s}\g_s(t)} \right) _{{\mathcal{H}}_\w(\g_s,t,u_s) } dt \right) _{t'=0}\\
&\quad \quad\quad \quad =\frac{d}{d t'} u_s g_{\g_s}\exp\left (t'\int_{0}^{1} \W \left (\t{\frac{\d}{\d t}\g_s(t)},\t{\frac{\d}{\d s}\g_s(t)} \right)   _{{\mathcal{H}}_\w(\g_s,t,u_s) } dt\right) _{t'=0},
\end{align*}
which, given the fact that the right action of $G$ on $P$ is free, is equivalent to:
\begin{Lemma}\label{EEE}
\begin{equation}\label{C}
\frac{d}{d s} g_{\g_s}=g_{\g_s} \int_{0}^{1} \W \left (\t{\frac{\d}{\d t}\g_s(t)},\t{\frac{\d}{\d s}\g_s(t)} \right) _{{\mathcal{H}}_\w(\g_s,t,u_s) } dt.
\end{equation}
\end{Lemma}
\noindent This fact will be of major importance later.

\subsection{Crossed modules}

All Lie groups and Lie algebras are taken to be finite-dimensional.

\begin{Definition}[Lie crossed module]\label{LCM}
 A crossed module ${\Gc= ( \d\colon E \to  G,\tr)}$ is given by a group morphism $\d\colon E \to G$ together with a  left action $\tr$ of $G$ on $E$ by automorphisms, {such that}:

\begin{enumerate}

 \item $\d(X \tr e)=X \d(e)X^{-1}; \forall X \in G, \forall e \in E,$

  \item $\d(e) \tr f=efe^{-1};\forall e,f  \in E.$

\end{enumerate}
If both $G$ and $E$ are Lie groups, $\d\colon E \to G$ is a smooth morphism,  and the left action of $G$ on $E$ is smooth then $\Gc$ will be called a Lie crossed module.

\end{Definition}

A morphism $\Gc \to \Gc'$ between the crossed modules ${\Gc= ( \d\colon E \to  G,\tr)}$  and {$\Gc'=(\d'\colon E' \to G',\tr')$} is given by a pair of maps $\f\colon G \to G'$ and  $\psi\colon E \to E'$ making the diagram:

$$ \begin{CD}
  E @>\d>> G \\
 @V \psi VV   @VV \f V \\
  E' @>\d'>> G' \\
   \end{CD}
$$
commutative. In addition we must have $\psi(X \tr e)=\f(X) \tr' \psi(e)$ for each $ e \in E$ and each $ X \in G$.

Given a Lie crossed module ${\Gc= ( \d\colon E \to  G,\tr)}$, then the induced Lie algebra map $\d\colon \le \to \lg$, together with the derived action of $\lg$ on $\le$ (also denoted by $\tr$) is a differential crossed module, in the sense of the following definition - see \cite{BS,B,BC}.

\begin{Definition}[Differential crossed module] A differential crossed module, say ${\mathfrak{G}=(\d \colon \le \to  \lg,\tr )}$, is given by a Lie algebra morphism $\d\colon \le \to \lg$ together with a left action of $\lg$ on the underlying vector space of  $\le$, {such that}:
\begin{enumerate}
 \item For any $X \in \lg$ the map $e \in \le \mapsto X \tr e \in \le$ is a derivation of $\le$, in other words $$X \tr [e,f]=[X \tr e,f]+[e, X \tr f];\forall X \in \lg, \forall e ,f\in \le.$$
\item The map $\lg \to \Der(\le)$ from $\lg$ into the derivation algebra of $\le$ induced by the action of $\lg$ on $\le$ is a Lie algebra morphism. In other words:
$$[X,Y] \tr e=X \tr (Y \tr e)-Y \tr(X \tr e); \forall X,Y \in \lg, \forall e\in \le.$$
\item $\d( X \tr e)= [X,\d(e)];\forall X \in \lg, \forall e\in \le.$
\item $\d(e) \tr f=[e,f];  \forall e,f\in \le.$
\end{enumerate}
\end{Definition}
\noindent Note that  the map $(X,e) \in \lg \times \le \mapsto X \tr e \in \le$  is necessarily bilinear.

Therefore, given a differential crossed module, ${\mathfrak{G}=(\d \colon \le \to  \lg,\tr )}$ there exists a unique crossed module of simply connected Lie groups ${\Gc=(\d \colon E \to  G,\tr )}$ whose differential form is $\mathfrak{G}$, up to isomorphism. The proof of this result is standard Lie theory, together with the lift of the Lie algebra action to a Lie group action, which can be found in \cite{K}, Thm 1.102.

\begin{Example}
Let $G$ be a Lie group with a left action $\tr$  on an abelian group $V$ by automorphisms, for example take $V$ to be any representation of $G$ on a vector space. {If we put $\d=1_G$ then 
$( \d \colon V \to G, \tr)$}
is a Lie crossed module. Its differential form is given by the derived action of $\lg$ on $V$ and the zero map $V \to \lg$.
\end{Example}

\begin{Example}
{Let $G$ be any Lie group. Let $\Ad$ denote the adjoint action of $G$ on $G$. Then ${(\id \colon G \to G, \Ad)}$ is a Lie crossed module.}
\end{Example}

\begin{Example}\label{w}
{Let $G$ and $E$ be Lie groups, and let $\d \colon E \to G$ be a surjective map such that $\ker(\d)$ is central in $E$. Define a left action of $G$ on $E$ as $g \tr e=g_0 e g_0^{-1}$, where $g_0\in E$ is such that $\d(g_0)=g\in G$, and $e \in E$. Then ${\Gc=(\d \colon E \to  G,\tr )}$ is a Lie crossed module.}
\end{Example}

\begin{Example}\label{r}
{In the previous example, we can take the universal covering map ${\rm SU}(2) \to {\rm SO}(3)$, whose kernel is $\{\pm 1\}$. The differential form of this crossed module is given by the adjoint action of ${\frak su}(2)$ on ${\frak su}(2)$ and the identity map.}
\end{Example}
\begin{Example}\label{cv}
{More generally, take the universal cover {$\d\colon E \to G$} of the Lie group $G$. Then the construction in example \ref{w} can be applied since the kernel of $\d$ is necessarily central in $E$.}
\end{Example}

\begin{Example}\label{Heis}
{Let $H$ be the group of upper triangular matrices in $\R$ which only have 1's in the main diagonal.  Define a Lie group map $\d\colon H  \to \R^2$ as
$$ \d \begin{pmatrix}1 & a & b \\ 0 &1 & c \\ 0 &0 &1 \end{pmatrix}=(a,c),$$
where $a,b,c \in \R$. Then the kernel of $\d\colon H\to \R^2$ is central in $H$, and therefore the construction in example \ref{w} can be applied.}
\end{Example}

\begin{Example}\label{q}
Let $E$ be any Lie group. Let $\Aut(E)$ be the group of {automorphisms} of $E$. Then $\Aut(E)$ is a Lie group, in fact it is a Lie subgroup of ${\rm GL}(\le)$ if $E$ is simply connected. Consider the map $\Ad\colon E \to \Aut(E)$ that sends $ e\in E$ to the automorphism $\Ad(e) \colon E \to E$. The group $\Aut(E)$ acts on $E$ as $\f \tr e=\f(e)$ where $\f\in \Aut(E)$ and $e \in E$. Then {$({\Ad}\colon E \to \Aut(E),\tr)$ is a Lie crossed module.}
\end{Example}
{The previous construction yields several  examples of Lie crossed modules. For instance, example \ref{r} is obtained this way since ${\rm Aut}({\rm SU(2)})={\rm SO(3)}$. To get new examples, however,  we need to consider non-semisimple Lie groups. }
\begin{Example}
{In example \ref{q}, take $E=H$, the group of upper triangular $3\times 3$ matrices in $\R$, which only have 1's in the main diagonal. Then the crossed module ${(\Ad\colon H \to  {{\Aut}}(H),\tr)}$ is highly non-trivial. In fact, $\Ad(E)$ is isomorphic to $\R^2$, whereas the cokernel ${\rm Out}(E)$ of $\Ad$ is isomorphic to ${\rm GL}(2,\R)$, see \cite{KH}. The kernel of $\Ad$ is a central subgroup of $H$ isomorphic to $\R$; see example \ref{Heis}.}
\end{Example}

\subsubsection{Crossed modules and categorical groups}\label{CCG}

\begin{Definition}[Categorical groups and Lie categorical groups]
A categorical group $\mathcal{C}$ is  a groupoid provided with a strict monoidal structure $\C \times \C \ra{\tn} \C$, as well as a ``group inversion'' functor ${^{-1}}\colon C  \to C$  such that both the set of objects $C_0$ of $\C$ and the set of morphisms $C_1$ of $\C$ are  groups under the tensor product, and both the source and target maps $\sigma,\tau \colon C_1 \to C_0$ are group morphisms. A Lie categorical group is defined analogously by replacing groups by Lie groups and group morphisms by  Lie group morphisms.
\end{Definition}

See \cite{BM,BL} for more details.

\begin{Definition}[Morphism of categorical groups]
A morphism of categorical groups is given by a strict monoidal functor; see \cite{ML}.
\end{Definition}

 It is well-known that the category of crossed modules and the category of categorical  groups are equivalent; see for example \cite{BL,BM,BHS}. Let us explain how to define a categorical group $\C(\Gc)$ from a crossed module $\Gc$.  This construction is an old one.

Let ${\Gc=(\d \colon E \to  G,\tr )}$ be a  crossed module.  
 The set of objects $C_0$ of $\C(\Gc)$ is given by all elements of $G$.  The set of morphisms $C_1$ of $\C(\Gc)$ is given by the set of all pairs $(X,e)$ where $X \in G$ and $e \in E$. The source and the target of $(X,e) \in C_1$ are given by $\sigma(X,e)=X$ and $\tau(X,e)=\d(e)^{-1}X$, respectively.
 In other words, a morphism in $C(\Gc)$  ``looks like''
$X \ra{(X,e)} \d(e)^{-1} X$, which we will sometimes abbreviate as $X \ra{e} \d(e)^{-1} X$. Given $X \in G$ and $e,f \in E$ the composition
$$X \ra{e} \d(e)^{-1}X  \ra{ f}   \d(f)^{-1}\d(e)^{-1} X$$ is $$X \ra{ef}  \d(ef)^{-1}X.$$   

The tensor product has the form:

\begin{equation}\label{Catmon}
\begin{CD}
\begin{CD}  \d(e)^{-1}X \\ @AA eA \\ X\\\end{CD} \tn \begin{CD}  \d(f)^{-1}Y \\ @AAfA \\ Y \end{CD}\quad  =\quad\begin{CD}\d(e)^{-1}X\d(f)^{-1}Y \\ @AA (X\tr f)e A \\ XY  \end{CD}
\end{CD},
\end{equation}
where $X,Y \in G$ and $e,f \in E$.
Therefore, the set of morphisms of $\C(\Gc)$ is a group under the tensor product. This group is isomorphic to the semidirect product $G \ltimes E$ of $G$ and $E$.

From the definition of a crossed module, it is easy to see that we have indeed defined a strict  {monoidal category which, furthermore,} is a categorical group.

We will recall in \ref{TCM} how to define a crossed module from a categorical group.


\subsection{The thin fundamental categorical group of a manifold}\label{thinfund}

Let $M$ be a smooth manifold.
\subsubsection{1-Tracks} \label{1-Tracks}

\begin{Definition}[1-path]
A 1-path is given by a smooth map $\gamma\colon [0,1] \to M$ such that there exists an $\epsilon >0$ such that $\g$ is constant in $[0,\e] \cup [1-\e,1]$; in the terminology of \cite{CP}, this can be abbreviated by saying that each end point of $\g$ has a sitting instant. Given a 1-path $\g$, define the source and target of $\g$ as $\sigma(\g)=\g(0)$ and $\tau(\g)=\g(1)$, respectively. 

\end{Definition}

Given two 1-paths $\g$ and $\f$ with $\tau(\g)=\sigma(\f)$, their concatenation $\g\f$ is  defined in the usual way:
$$(\g\f)(t)=\left \{ \begin{CD} \g(2t), \textrm{ if }t \in [0,1/2] \\ \f(2t-1), \textrm{ if } t \in [1/2,1]\end{CD} \right.$$ 
Note that the concatenation of two 1-paths is also a 1-path.  The fact that any 1-path has sitting instants at its end points needs to be used to prove this.

\begin{Definition}[2-paths]

A 2-path $\G$ is given by a smooth map $\G\colon [0,1]^2 \to M$ such that 
there exists an $\epsilon >0$ for which:
\begin{enumerate}

\item $\G(t,s)=\G(0,0)$ if $0 \leq t \leq \epsilon$ and $s \in [0,1]$,

\item $\G(t,s)=\G(1,0)$ if  $1-\epsilon \leq t \leq 1$ and $s \in [0,1]$,

\item $\G(t,s)=\G(t,0)$ if $0 \leq s \leq \epsilon$ and $t \in [0,1]$,

\item $\G(t,s)=\G(t,1)$ if  $1-\epsilon \leq s \leq 1$ and $t \in [0,1]$.

\end{enumerate}

\end{Definition}

{Given a 2-path $\Gamma$, define the following 1-paths:}
\begin{align*}
\d_l(\G)(s)&=\G(0,s), s \in [0,1],
&\d_r(\G)(s)&=\G(1,s),  s \in [0,1],\\
\d_d(\G)(t)&=\G(t,0), t \in [0,1],
&\d_u(\G)(t)&=\G(t,1), t \in [0,1].
\end{align*}
Note that {the paths $\d_l(\G)$ and $\d_r(\G)$} are necessarily constant, thus each of them can be identified with a  point of $M$. In addition, given a 2-path $\G$ then  the 1-paths $\d_u(\G)$ and $\d_d(\G)$ have the same initial and end points. {If $x$ is a point of $M$ the constant 1- and 2-paths with value $x$ are both denoted by $x$.}

If $\G$ and $\G'$ are 2-paths such that $\d_r(\G)=\d_l(\G')$ their horizontal concatenation $\G \circ_h \G'$ is defined in the obvious way, in other words:

$$  \big (\G \circ_h \G' \big )(t,s)=\left \{ \begin{CD} \G(2t,s), \textrm{ if }t \in [0,1/2] \textrm{ and }  s \in [0,1] \\ \G'(2t-1,s), \textrm{ if } t \in [1/2,1] \textrm{ and }  s \in [0,1] \end{CD} \right. $$ 
Similarly, if $\d_u(\G)=\d_d(\G')$ we can define a vertical concatenation $\G \circ_v \G'$ as:
$$\big(\G \circ_v \G' \big)(t,s)=\left \{ \begin{CD} \G(t,2s), \textrm{ if }s \in [0,1/2]  \textrm{ and }  t \in [0,1]\\ \G'(t,2s-1), \textrm{ if } s \in [1/2,1]  \textrm{ and }  t \in [0,1] \end{CD} \right.$$

We will also represent the  $2$-paths $\G$ of $M$ in the following  suggestive way:
$$\G=\begin{CD} \d_u(\G )\\ @AA\G A \\ \d_d(\G) \end{CD}$$
In this notation, the vertical concatenation of two $2$-paths $\G$ and $\G'$ with $\d_d(\G')=\d_u(\G)$ can be represented as:

$$
     \G \circ_v \G'\quad=\quad \begin{CD} \d_u(\G\circ_v \G' )\\ @AA\G\circ_v \G' A \\ \d_d(\G\circ_v \G' ) \end{CD}\quad =\quad \begin{CD} \d_u(\G' )\\ @AA\G' A\\  \d_d(\G')=\d_u(\G)\\ @AA\G A\\ \d_d(\G) \end{CD}
$$

\begin{Definition}
Two 1-paths $\f$ and $\g$ are said to be rank-1 homotopic  (and we write $\f\cong_1 \g$) if there exists a 2-path $\G$ such that:
\begin{enumerate}
 \item $\d_u(\G)=\g$ and $\d_d(\G)=\f$.
  \item $\Rank({{\EuScript{D}}}_{v} \G) \leq 1, \forall v \in [0,1]^2.$
\end{enumerate}
{Here ${\EuScript{D}}$ denotes derivative.}
\end{Definition}
In particular, if $\g$ and $\f$ are rank-1 homotopic, then they have the same initial and end-points.  Note also that rank-1 homotopy is an equivalence relation. Given a 1-path $\g$, the equivalence class to which it belongs is denoted by $[\g]$. Rank-1 homotopy is one of a number of notions of ``thin'' equivalence between paths or loops, and was introduced in 
\cite{CP}, following a suggestion by A. Machado.

We denote the set of $1$-paths of $M$ by $S_1(M)$. The quotient of $S_1(M)$ by the relation of thin homotopy is denoted by $\S_1(M)$. We call the elements of $\S_1(M)$  1-tracks. 
 The {concatenation} of $1$-tracks together with the source and target maps $\sigma, \tau \colon \S_1(M) \to M$, defines a groupoid $\Sc_1(M)$ whose set of morphisms is  $\S_1(M)$ and whose set of objects is $M$.

\begin{Definition}

Let $* \in M$ be a base point. The group $\pi_1^1(M,*)$ is defined as being the set of  1-tracks $[\g]\in \S_1(M)$ starting and ending at $*$, with the group operation being the concatenation of paths. 

\end{Definition}
It is easy to see directly that $\pi_1^1(M,*)$, a thin analogue of $\pi_1(M)$, is indeed a group. See \cite{CP,MP} for details. 

\subsubsection{2-Tracks}\label{2-Tracks}

\begin{Definition}\label{thin2}

Two 2-paths $\G$ and $\G'$ are said to be rank-2 homotopic  (and we write $\G\cong_2 \G'$) if there exists a smooth map ${J}\colon [0,1]^2 \times [0,1] \to M$ such that:

\begin{enumerate}

\item {$J(t,s,0)=\G(t,s),\, J(t,s,1)=\G'(t,s)$ for $s,t\in [0,1]$.} 
 \item ${J}$ is constant over $\{0\} \times [0,1]^2$ and  $\{1\} \times [0,1]^2$.

 \item Over $[0,1] \times \{0\} \times [0,1]$ the map ${J} \colon [0,1]^2 \times [0,1] \to M$ restricts to a rank-1 homotopy $\d_d(\G) \to \d_d(\G')$.  

 \item ${J}$ restricts to a rank-1  homotopy $\d_u(\G) \to \d_u(\G')$ over $[0,1] \times \{1\} \times [0,1]$.

\item

 There exists an $\epsilon >0$ such that $J(t,s,x)={J}(t,s,0)$ if $x\leq \epsilon$ and $s,t \in [0,1]$, and analogously for all the other faces of $[0,1]^3$.  We will denote this condition as saying that $J$ has a product structure close to the boundary of $[0,1]^3$.

 \item  $\Rank ( {{\EuScript{D}}}_{v} {J}) \leq 2$ for any $v \in {[0,1]}^3$. 

\end{enumerate}

\end{Definition}

Note that {rank-2 homotopy} is an equivalence relation. To prove {transitivity} we need to use the penultimate condition of the previous definition. 
We denote by $S_2(M)$ the set of all 2-paths of $M$. The quotient of $S_2(M)$ by the relation of rank-2 homotopy is denoted by $\S_2(M)$. We call the elements of $\S_2(M)$ 2-tracks.
 If $\G \in S_2(M)$, we denote the equivalence class in $\S_2(M)$ to which $\G$ belongs by $[\G]$.

\subsubsection{Horizontal and vertical compositions of 2-tracks}\label{HVC}

Note that the horizontal composition of 2-paths descends immediately to a horizontal composition in $\S_2(M)$.  Suppose that $\G$ and $\G'$ are 2-paths with {$\d_u(\G)\cong_1 \d_d(\G')$.} Choose a rank-1 homotopy ${J}$ connecting $\d_u(\G)$ and $\d_d(\G')$. Then $[\G]\circ_v[\G']$ is defined as {$[(\G \circ_v J) \circ_v \G'].$} The fact that this composition is well defined in $\S_2(M)$ follows immediately from  lemma \ref{MAIN} in the appendix.

Therefore both the vertical and horizontal compositions of 2-paths descend to $\S_2(M)$.  These compositions are obviously associative, and admit units and inverses. Since the interchange law is trivially verified the following theorem holds:
\begin{Theorem}Let $M$ be a smooth manifold.
The horizontal and vertical {compositions}  in $\S_2(M)$ together with the boundary maps $\d_u,\d_d,\d_l,\d_r\colon \S_2(M) \to \S_1(M)$ define a 2-groupoid $\Sc_2(M)$,  whose set of objects is given by all points of $M$, the set of 1-morphisms by the set $\S_1(M)$ of 1-tracks on $M$, and given two 1-tracks $[\g]$ and $[\g']$, the set of 2-morphisms $[\g] \to [\g']$ is given by all 2-tracks $[\G]$ with $\d_d([\G])=[\g]$ and $\d_u([\G])=[\g']$.
\end{Theorem}
The definition of a 2-groupoid can be found in \cite{HKK}.
This construction should be compared with \cite{HKK,BHKP}, where the thin strict 2-groupoid of a Hausdorff space was defined, using a different notion of thin equivalence (factoring through a graph) - see also \cite{BH} for the construction of the fundamental double groupoid of a triple of spaces. For analogous non-strict constructions see \cite{M,BS,MP}.

\subsubsection{Definition of the thin fundamental categorical group of a smooth manifold}\label{TCM}

Since $\Sc_2(M)$ is a 2-category, given a base point $*$ of $M$ (in other words an object of $\Sc_2(M)$) we can define a monoidal category $\Pc_2(M,*)$. The objects of $\Pc_2(M,*)$ are given by all  1-tracks of $M$ starting and ending at $*$, and therefore by all elements of $\pi_1^1(M,*)$. On the other hand, the set $\P_2(M,*)$ of {morphisms} of $\Pc_2(M,*)$ is given by the set of all 2-tracks $[\G] \in \S_2(M)$ such that $\d_l([\G])=\d_r([\G])=*$. 
Therefore,  given $[\g],[\g'] \in \pi_1^1(M,*)$, the set of morphisms $[\g] \to [\g']$ is given by all 2-tracks $[\G]$ connecting $[\g]$ and $ [\g']$, in other words with $\d_{{d}}([\G])=[\g]$ and $\d_{{u}}([\G])=[\g']$.  The composition in   $\P_2(M,*)$ is given by the vertical composition in $\S_2(M)$. On the other hand the horizontal composition of 2-tracks will give the tensor product in $\Pc_2(M,*)$.

 Given that $\Sc_2(M)$ is a 2-groupoid, the monoidal category $\Pc_2(M,*)$ is a categorical group. We therefore call $\Pc_2(M,*)$  the {\it thin fundamental categorical group of $M$.}

\subsubsection{The thin fundamental crossed module of a smooth manifold}\label{fcm}
Since   $\Pc_2(M,*)$   is a categorical group, it  defines a crossed module $\Pi_2^{2,1}(M,*)$, which we call the {\it thin fundamental crossed module of $M$.} Let us explain how to construct $\Pi_2^{2,1}(M,*)$ directly.
Define a group $\pi_2^{2,1}(M,*)$ as being given by 
$$\pi_2^{2,1}(M,*)=\left \{[\G] \in \P_2(M,*) \colon \d_d({[\G]})=[*]\right \}$$
with the  product law being the horizontal composition. Note that the map $\d \doteq \d_u \colon \pi_2^{2,1}(M,*)\to \pi_1^1(M,*)$ is therefore a group morphism. In addition $\pi_1^1(M,*)$ acts on $\pi_2^{2,1}(M,*)$ on the left as follows:
$$ [\g] \tr [\G]=\id_{[\g]} \circ_h [\G] \circ_h\id_{[\g^{-1}]},$$
where $[\G]\in {\pi_2^{2,1}(M,*)}$ and $[\g] \in \pi_1^1(M,*). $
Here if $[\g]$ is a 1-track then  $\id_{[\g]}$ denotes the 2-track such that $\id_{[\g]}(s,t)=\g({t})$, where $s,t \in [0,1]$.

The following theorem follows from the well-known equivalence between the category of  categorical groups and  the category  of crossed modules which was referred to in \ref{CCG}; see \cite{BM,BL,BMo,P}. For a graphical proof see \cite[6.2]{BHS}.
\begin{Theorem}
The boundary map $\d \colon \pi_2^{2,1}(M,*)\to \pi_1^1(M,*)$, together with the left action of  $\pi_1^1(M,*)$  on $\pi_2^{2,1}(M,*)$, defines a crossed module, which we  denote by $\Pi_2^{2,1}(M,*)$. 
\end{Theorem}

In \ref{CCG} we showed how to construct a categorical group from any crossed module.  The way we defined $\Pi_2^{2,1}(M,*)$ from the categorical group $\Pc_2(M,*)$ shows how to go in the reverse direction.

\subsubsection{Two exact sequences}\label{2es}
This subsection is not necessary to understand the rest of the article.

\begin{Definition}
Let $M$ be a smooth manifold. Let $* \in M$. The group $\pi_2^2(M,*)$ is given by all 2-paths $\G \in S_2(M)$ such that $\G(\d([0,1]^2))=\{*\}$, up to the equivalence relation given by $\G \cong \G'$ if there exists a rank-2 homotopy ${J\colon [0,1]^2\times [0,1] \to M}$, connecting $\G$ and $\G'$, such that, in addition, $J\left ((\d ([0,1]^2) \times I\right)=\{*\}$. The group law is taken to be the horizontal composition.
\end{Definition}
For more details see \cite{MP}.

\begin{Theorem}
Let $M$ be a smooth manifold {with a base point $*$}. The crossed module $\Pi_2^{2,1}(M,*)$  can be embedded into the exact sequence:
$$\{0\} \to \pi_2^2(M,*) \ra{i} \pi_2^{2,1}(M,*) \ra{\d} \pi_1^1(M,*) \ra{p} \pi_1(M,*)\to \{1\}.$$
\end{Theorem}
\begin{Proof}
The only non-trivial part is the proof that the natural  map $i\colon \pi_2^2(M,*)  \to \pi_2^{2,1}(M,*)$ is injective. This follows from Lemma \ref{MAIN}.
\end{Proof}

\begin{Definition}
 The thin $k$-{invariant} of $M$ is defined as being the cohomology class {$k_t(M,*) \in H^3\left (\pi_1(M,*),\pi_2^2(M,*)\right)$} determined by the exact sequence above. See \cite{EML,Br} for the construction of these group cohomology classes.
\end{Definition}

We can define a 2-category $\Qc_2(M)$, with objects given by the points of $M$ and  morphisms given by the arrows of the groupoid $\S_1(M)$. The set of 2-morphisms $\Q_2(M)$ of $\Qc_2(M)$ is given by  the set of $2$-paths of $M$ up to homotopy. Here the notion of  homotopic 2-paths  is defined in the same way as the notion of rank-2 homotopic 2-paths (Definition \ref{thin2}), but omitting the final condition on the rank of {the derivative of} the map defining the homotopy; see the remarks after definition \ref{K}. If $\G$ and $\G'$ are 2-paths, we write $\G \cong_\infty \G'$ to denote that they are homotopic. Applying the construction that defined the crossed module $\Pi_2^{2,1}(M,*)$ to $\Qc_2(M,*)$ yields a crossed module $$\Pi_2^{\infty,1}(M,*)=\left ( \pi_2^{\infty,1}(M,*) \ra {\d} \pi_1^1(M,*), \tr \right).$$
Here the group $\pi_2^{\infty,1}(M,*)$ is given by 
$$\pi_2^{\infty,1}(M,*)=\left \{[\G] \in \Qc_2(M,*) \colon \d_d({[\G]})=[*]\right \}$$
with the  product law being the horizontal composition. Therefore the  map $\d \doteq \d_u \colon {\pi_2^{\infty,1}(M,*)}\to \pi_1^1(M,*)$ is  a group morphism.

As before we have:

\begin{Theorem}
The crossed module $\Pi_2^{\infty,1}(M,*)$  can be embedded into the exact sequence:
$$\{0\} \to \pi_2(M,*) \ra{i} \pi_2^{\infty,1}(M,*) \ra{\d} \pi_1^1(M,*) \ra{p} \pi_1(M,*)\to \{1\}.$$
\end{Theorem}
{It is easy to show that the  cohomology class $k(M,*) \in  H^3\left (\pi_1(M,*),\pi_2(M,*)\right)$ given by this exact sequence coincides with the $k$-{invariant} of $M$. This follows from the explicit construction of it in \cite{EML}, and also from  the no-free-lunch principle.}

\section{Two-dimensional holonomies}

\subsection{Categorical connections}
\begin{Definition}
Let $G$ be a Lie group. Let $M$ be a smooth manifold.  Let $P \to M$ be a smooth principal  bundle with structure group $G$. Denote the right action of $G$ on $P$ as $g \in G\mapsto {R_g} \in \diff(M)$. Let also ${\Gc=(\d\colon E \to G,\tr)}$  
be a Lie crossed module, where $\tr$ is a Lie group left action of $G$ on $E$ by automorphisms. Let also  ${\mathfrak{G}=(\d \colon \le \to  \lg,\tr )}$ be the associated differential crossed module. A $\Gc$-categorical connection on $P$ is a pair $(\w,m)$, where $\w$ is a connection 1-form on $P$, {i.e.} $\w\in \A^1(P,\lg)$ is a $1$-form on $P$ with values in $\lg$ such that: 
\begin{enumerate}
\item $R_g^*(\w)=g^{-1}\w g , \forall g \in G,$
\item $\w(A^\#)=A,\forall A \in \lg$;
\end{enumerate}
and $m\in \A^2(P,\le)$ is a 2-form on $P$ with values in $\le$, the Lie algebra of $E$, such that:\begin{enumerate}
\item The 2-form $m$ is $G$-{equivariant}. In other words we have $R_g^*(m)=g^{-1} \tr m$ for each $g \in G$. 
\item The 2-form $m$ is horizontal, in other words:
 $$m(X,Y)=m(X^H,Y^H), \textrm{ for each } X,Y \in \X(P).$$ 
In particular $m(X_u,Y_u)=0$ if either of the vectors $ X_u,Y_u \in T_u P$ is  vertical, where $u \in P$. Here the map $X \in \X(P) \mapsto X^H \in \X(P)$ denotes the horizontal projection of vector fields on $P$ with respect to the connection 1-form $\w$.
\end{enumerate}
Finally we require that 
\begin{equation}
\d(m)=\W
\end{equation} 
where $\W\in \A^2(P,\lg)$ is the curvature 2-form of $\w$.
\end{Definition}
\noindent{This last condition} is of course equivalent to the vanishing of the fake curvature of \cite{BS,BrMe}; see \ref{lform}.
\begin{Example}{
Let $P$ be any principal $G$-bundle with a connection 1-form $\w$. Let $\W$ be the curvature $2$-form of $\w$. If $\Gc=(\id \colon G \to G,\tr)$, where $\tr$ is the adjoint action of $G$ on $G$, then $(\w,\W)$ is a $\Gc$-categorical connection on $P$. More generally we can take $\Gc$ to be given by the universal covering map $E \to G$ of $G$; {see example \ref{cv}.} }\end{Example}

\subsection{The {2-curvature 3-form} of a $\Gc$-categorical connection}
Recall the definition of the curvature 2-form $\W$ of a connection 1-form $\w$ as the covariant exterior derivative of the 1-form $\w$.

\begin{Definition}[{2-Curvature}]
Let ${\Gc=(\d\colon E \to G,\tr)}$  be a Lie crossed module, and let $P \to  M$ be  a smooth principal $G$-bundle. 
The {2-curvature} 3-form of a $\Gc$-categorical connection $(\w,m)$  on $P$ is defined as 
$$\M=d m\circ(H \times H \times H).$$
In other words:
$$\M(A,B,C)=dm(A^H,B^H,C^H),$$
where $A,B,C \in \X(P)$ are smooth vector fields on $P$. 
\end{Definition}

\subsubsection{Algebraic preliminaries}
Let $M$ be a  smooth manifold. Consider a differential crossed module ${\mathfrak{G}=(\d\colon \le \to \lg,\tr)}$. In particular the map $(X,e) \in \lg \times \le  \mapsto X \tr e \in \le$ is bilinear. 

Let $a \in \A^n(M,\lg )$  and $b \in \A^m(M,\le)$ be $\lg$- and $\le$-valued (respectively) differential forms on  $M$. We define $a \tn^\tr b  $ as being the  $\le$-valued covariant tensor field on $M$ such that
$$(a \tn^\tr b)(A_1,\ldots, A_n,B_1,\ldots B_m)= a (A_1,\ldots, A_n)\tr b(B_1,\ldots, B_m);A_i,B_j \in \X(M).$$
{We also define an alternating tensor field  $a \wedge^\tr b \in \A^{n+m}(M,\le)$, being given by }
$${ a \wedge^\tr b=\frac{(n+m)!}{n! m!}{\rm Alt}(a  \tn^\tr b).}$$
Here ${\rm Alt}$ denotes the natural projection from the vector space of $\le$-valued covariant tensor fields on $M$ onto the vector space of $\le$-valued differential forms on $M$.

For example, if  $a \in \A^1(M,\lg)$ and $b \in \A^2(M,\le)$, then $a \wedge^\tr b$ satisfies:
\begin{equation}\label{BB}
(a \wedge^\tr b)(X,Y,Z)=a(X) \tr b(Y,Z)+a(Y) \tr b(Z,X)+a(Z) \tr b(X,Y),
\end{equation}
where $X,Y,Z \in \X(M).$  

Note that the formula:
$$d (a  \wedge^\tr b)=d a  \wedge^\tr b  +(-1)^n a  \wedge^\tr d b,$$
where $a \in \A^n(M,\lg)$ and $b \in \A^m(M,\le)$, 
holds trivially.

As as example of this notation, let $P\to M$ be a principal $G$-bundle over the manifold $M$. Let $\w \in \A^1(P,\lg)$ be a connection 1-form on $M$. Then the structure equation  for the curvature form $\W=D \w$ of $\w$ can be written as:
$$\W=d\w+{\frac{1}{2}}\w \wedge^{\ad} \w,$$
where $\ad$ denotes the adjoint action of $\lg$ on $\lg$.
Furthermore the Bianchi identity can be written as:
$$d \W+\w\wedge^{\ad} \W=0;$$
see \ref{curv}.

\subsubsection{{The 2-structure equation}} \label{se}
The following equation is an analogue of Cartan's structure equation. It will be of prime importance later.
\begin{Proposition}[{2-Structure Equation}]
Let ${\Gc= ( \d\colon E \to  G,\tr)}$ be a Lie crossed module. Let $P$ be a principal $G$-bundle over the manifold $M$. Let also $(\w,m)$ be a categorical $\Gc$-connection on $P$. We have:
$\M=dm {+}\w\wedge^\tr m,$ where {$\M=D m$} is the {2-curvature} $3$-form of $(\w,m)$. {In particular {the} 2-curvature 3-form $\M$ is $G$-{equivariant}, in other words: $R_g^*(\M)=g^{-1} \tr \M,$
for each $g \in G$.}
\end{Proposition}

The proof of the 2-structure equation follows directly from the following lemma.
\begin{Lemma}\label{ext}
Let $a$ be a $G$-{equivariant} horizontal {$(n-1)$-form} in $P$.  Then $Da=da+\w\wedge^\tr a$.
\end{Lemma}
\begin{Proof}
{Let $X^1,\ldots,X^n \in T_uP$, extended to a neighbourhood of $u\in P$. If all vectors are horizontal it follows that $Da(X^1,\ldots,X^n)=(da+\w\wedge^\tr a)(X^1,\ldots,X^n)$, {since $\w$ is vertical.} Suppose that one of the vectors, which we can suppose {to be}  $X^1$, is vertical. We need to prove that $(da+\w\wedge^\tr a)(X^1,\ldots,X^n)=0$.  Choose $A \in \lg$ such that $A^\#_u=X^1_u.$ By using {Cartan's} magic formula together with the fact that $a$ is horizontal {it} follows that {$da (A^\#,X^2,\ldots,X^n)=(\L_{A^\#}a) (X^2,\ldots X^n)$}, where $\L$ denotes Lie derivative. Since $a$ is $G$-{equivariant} {it} follows that $\L_{A^\#}a=-A \tr a $. {However, $(\w\wedge^\tr a)(A^\#,X^2,\ldots,X^n)=\w(A^\#)\tr a(X^2,\ldots,X^n)$, since $a$ is horizontal. Note $\w(A^\#)=A$.}}
\end{Proof}

\subsubsection{{The 2-Bianchi identity}}
Fix a Lie crossed module  ${\Gc= ( \d\colon E \to  G,\tr)}$.
Let $P$ be a smooth principal $G$-bundle over $M$ with a connection $\w$. Choose a $\Gc$-categorical connection $(\w,m)$ on $P$.

\begin{Proposition}[{2-Bianchi identity}]
Let $\M\in \A^3(P,\le)$ be the {2-curvature} 3-form of $(\w,m)$. Then the exterior covariant derivative $D \M$ of $M$ vanishes. In other words:
$$d \M \circ (H \times H \times H \times H)=0,$$
which by lemma \ref{ext} is the same as:
$$d \M +\w \wedge^{\tr} \M=0.$$
\end{Proposition}

\begin{Proof}
{Let us  prove the second form of the 2-Bianchi identity. By the 2-structure equation it follows that:}
$$
d \M=dd m{+}d(\w \wedge ^\tr m)=d \w \wedge ^\tr m{-} \w \wedge ^\tr d m 
=d \w \wedge ^\tr m{-} \w \wedge ^\tr \M + \w\wedge^{\tr} \left (\w \wedge^{\tr} m\right).$$
Given that $\w\wedge^{\tr} \left (\w \wedge^{\tr} m\right)=\frac{1}{2}(\w\wedge^{\ad } \w )\wedge ^\tr m$ we can conclude:
$$d \M=(d \w+\frac{1}{2}\w\wedge^{\ad } \w ) \wedge ^\tr m{-} \w \wedge ^\tr \M 
=\W \wedge ^\tr m{-} \w \wedge ^\tr \M 
=\d(m) \wedge ^\tr m{-} \w \wedge ^\tr \M.$$
We now prove that $\d(m) \wedge ^\tr m=0$. Since ${\mathfrak{G}=(\d \colon \le \to  \lg,\tr )}$ is a differential crossed module we have, for any vector fields $X,Y,Z,W \in \X(P)$:$$
{\d(m)(X,Y) \tr m(Z,W)=\d\big (m(X,Y)\big ) \tr m(Z,W)
                     =[m(X,Y), m(Z,W)].}
$$
Let us see  that ${\rm Alt}[m,m]=0$, which implies that $\d(m) \wedge ^\tr m=0$. The permutation $r$ such that $r(X,Y,Z,W)=(Z,W,X,Y)$ is {even}. In particular ${\rm Alt}([m,m])= {\rm Alt}([m,m] \circ r)$. However, $[m,m] \circ r=-[m,m]$.
\end{Proof}

\subsubsection{{Local form of a categorical connection}}\label{lform}

Let ${\Gc=(\d\colon E \to G,\tr)}$ be some Lie crossed module, and let ${\mathfrak{G}=(\d \colon \le \to  \lg,\tr )}$ be the associated differential crossed module. Fix a smooth manifold $M$ and a smooth principal $G$-bundle $P \to M$, with a categorical $\Gc$-connection $(\w,m)$, where, as usual, $\w \in \A^1(P,\lg)$ and $m \in \A^2(P,\le)$.

The local form of the categorical $\Gc$-connection {$(\w,m)$} is similar to the local form of a connection. Let $\{U_\a\}$ be an open covering of the manifold $M$. We can suppose that each $U_\a$ is contractible and that, for each $\a$ and $\b$, the intersection $U_\a \cap U_\b$ is contractible; see \cite{BT}. Therefore, for each $\a$, the restriction $P_\a$ of $P$ to $U_\a$ is trivial, and, in particular, it  admits a section $f_\a\colon U_\a \to P_\a$. The local form of a categorical $\Gc$-connection is therefore given by the family of forms $\w_a=f_\a^*(\w)$ and $m_\a=f_\a^*(m)$. 

For each $\a$, by the structure equation, the local curvature form $\W_\a=f_\a^*(\W)$ takes the form $$\W_\a=d \w_\a+\frac{1}{2} \w_\a\wedge^{\ad}\w_\a.$$ Locally, the Bianchi identity reads:
$$d \W_\a+\w_\a\wedge^\ad \W_\a=0. $$
Analogously, from the 2-structure equation, locally in each $U_\a$, the 2-curvature 3-form  $\M_\a=f_\a^*(\M)$ takes the form: $$\M_\a=d m_\a+\w_\a\wedge^{\tr}  m_\a,$$
and the 2-Bianchi identity reads:
$$d \M_\a+\w_\a \wedge^{\tr} \M_\a=0.$$
{These local expressions correspond to the formulae in \cite{B} when the fake curvature vanishes.}

Let $\f_{\a,\b} \colon U_a \cap U_\b \to G$ be defined as $f_\a(x)\f_{\a,\b}(x)=f_\b(x)$, for each $x \in  U_\a \cap U_\b$. 
We have:
\begin{align}
\label{x} \w_\b(X)&=\f_{\a,\b}^{-1} \w(X) \f_{\a,\b}+\f_{\a,\b} ^{-1} d \f_{\a,\b},\\
\label{y} \W_\b(X,Y)&=\f_{\a,\b}^{-1}  \W (X,Y)\f_{\a,\b}, 
\end{align}
for each $X,Y \in \X(M)$. Since both $m$  and $\M$ are  $G$-{equivariant} and horizontal, it follows that for each $\a,\b$:
\begin{align}\label{z}
 m_\b(X,Y)&=\f_{\a,\b}^{-1} \tr  m_\a(X,Y),\\
 \M_\b(X,Y,Z)&=\f_{\a,\b}^{-1} \tr  \M_\a(X,Y,Z),
\end{align}
for each $X,Y,Z \in \X(M)$.

\begin{Lemma}
Let $M$ be a smooth manifold. Let also ${\Gc=(\d\colon E \to G,\tr)}$ be a Lie crossed module, and let ${\mathfrak{G}=(\d \colon \le \to  \lg,\tr )}$ be the associated differential crossed module.
Suppose $\pi\colon P \to M$ is the trivial bundle $M \times G$. Given  $m_0 \in \A^2(M,\le)$ and $\w_0 \in \A^1(M,\lg)$, such that {$\d(m_0)=\W_0=d \w_0+\frac{1}{2} \w_0 \wedge^{\ad} \w_0$},  there exists a  categorical $\Gc$-connection $(\w,m)$ on $P$ such that its local form is $(\w_0,m_0)$.  
\end{Lemma}
\begin{Proof}
Let $\theta$ be the canonical $\lg$-valued left-{invariant} 1-form on $G$; see \ref{wkl}. Then $\w=g^{-1}\pi^*(\w_0)g+\theta$ is a connection 1-form on $M \times G$, and its local form, considering the section $f\colon M \to P$ such that, for each $x \in M$, $f(x)=(x,1_G)$,  is $\w_0$.
The 2-form $\pi^*(m_0)\in \A^2(P,\le)$ is horizontal. However, it is not $G$-{equivariant}. Therefore we define $m\in \A^2(M \times G,\le)$ as being such that, given $X,Y \in T_{(x,g)}P$, we have:
$$m(X,Y)=g^{-1}\tr \pi^*(m_0)(X,Y),$$
where $g \in G$, $x \in M$ and $X,Y \in T_{(x,g)} M\times G$. 
This 2-form is obviously smooth and horizontal. Moreover, given $X,Y \in T_{(x,g)}P$ it follows that:
\begin{align*}
(R_h^*m)(X,Y)&=m(Xh,Yh)\\
             &=(gh)^{-1}\tr \pi^*(m_0)(Xh,Yh)\\
             &=h^{-1}g^{-1}\tr \pi^*(m_0)(X,Y)   \\
             &=h^{-1}\tr m(X,Y).
\end{align*}
Therefore $m$ is $G$-{equivariant}.  Let $\W$ be the curvature of $\w$. Then  $\W(X,Y)=g^{-1} \pi^*(\W_0)(X,Y)g$, where as before $X,Y \in T_{(x,g)}P$. In particular it follows that $\W=\d(m)$.
\end{Proof}

The following follows from the previous discussion.
\begin{Corollary}
Let $M$ be a smooth manifold. Let also ${\Gc=(\d\colon E \to G,\tr)}$ be a Lie crossed module, and let ${\mathfrak{G}=(\d \colon \le \to  \lg,\tr )}$ be the associated differential crossed module. Let $\{U_\a\}$ be an open cover of $M$, {such} that the restriction $P_\a$ of $P$ to $U_\a$ is the trivial bundle. Let also $f_\a \colon U_\a \to P_\a$ be local sections of $P$. Define $\f_{\a,\b} \colon U_\a \cap U_\b \to G$ as $f_\a(x)\f_{\a,\b}(x)=f_\b(x)$. Suppose that for each $\a$  we are given forms $\w_\a \in \A^1(U_\a,\lg)$ and $m_\a \in \A^2(U_\a,\le)$, satisfying $\d(m_\a)=\W_\a={{d\w_\a}+\frac{1}{2}\w_\a \wedge^{\ad} \w_\a}$. Moreover, suppose that for each $\a,\b$ with $U_\a \cap U_\b \neq \emptyset$ equations (\ref{x}) and (\ref{z}) are satisfied. There exists a unique categorical connection {$(\w,m)$} on $P$ such that the forms $\w_\a$ and $m_\a$ are its local form.
 \end{Corollary}

The trivial bundle $M \times G\to M$ always admits the trivial categorical connection, whose local form is given by the forms $\w_0=0$ and $m_0=0$. However, there exist several other examples for the trivial bundle. Let $A=\d(E)$, and let $\la$ be its Lie algebra. There exists a linear map $r\colon \la \to \le$ such that $\d \circ r=\id_\la$. Let $\w_0$ be any $\la$-valued 1-form on $M$. Let $\W_0$ be the curvature of $\w_0$ and define $m_0=r(\W_0)$. Then $(\w_0,m_0)$ is the local form of a categorical $\Gc$-connection on $M\times G$. 
In fact, let $B=\ker\d\subset E$. Therefore $B$ is closed under the action of $G$ on $E$, by the first condition of the definition of crossed modules; see  definition \ref{LCM}. Let $\lb$ be the Lie algebra of $B$. Given any $\lb$-valued 2-form $\rho_0$ on $M$ then $\rho_0+m_0$ is the local form of a categorical $\Gc$-connection on $M\times G$.

It is not a trivial problem whether there always exists a categorical
 $\Gc$-connection on any  principal $G$-bundle. This issue is very likely to depend on the topology of the  manifold.
The strategy of combining locally defined connection 1-forms $\w$  using a partition of unity cannot be carried through for the 2-forms $m$, since $\partial (m) = \W$, which is not {${\rm C}^\infty$-linear in $\w$. }

\subsection{Categorical holonomy}
\subsubsection{Definition of a holonomy.}
The following definition appears in \cite{CP}.
\begin{Definition}[Holonomy]
 {Let $M$ and $N$ be smooth manifolds. Suppose that  $M$ has a chosen base point $*$.  A map $F \colon \pi_1^1(M,*) \to N$ is said to be smooth if, for any  $2$-path $\G\colon [0,1]^2 \to M,$ with $\d_l(\G)=\d_r(\G)=*$, the map $$s \in [0,1] \mapsto F([\g_s])\in N$$ is smooth.  
Here $\G(t,s)=\g_s(t)$, where  $s,t \in [0,1]$.}

 {Let $G$ be a Lie group. A holonomy is, by definition, a smooth group morphism $\pi_1^1(M,*) \to G$.}
\end{Definition}

The following result appears in \cite{CP}.

\begin{Theorem}\label{E}
Let $M$ be a smooth manifold with a base point $*$. Let also $G$ be a Lie group. Let $P \to M$ be a smooth principal $G$-bundle over $M$ with a connection $\w$. Let $u \in P_*$, the fibre of $P$ at $*$. Then the parallel transport of  $\w$ determines {a holonomy} $F^1_{(\w,u)} \colon \pi_1^1(M,*) \to G$. 
\end{Theorem}
\begin{Proof}
Recall the notation of \ref{Hol}. Given a 1-path $\g \colon [0,1] \to M$, with $\g(0),\g(1)=*$, define  $F^1_{(\w,u)}(\g)$ {to be} the unique element of $G$ such that $u={\mathcal{H}}_\w(\g,1,u)F^1_{(\w,u)}(\g)$.  Then, if $\G$  is a 2-path, with $\d_l(\G)=\d_r(\G)=*$, it follows that the map 
$$s \in [0,1] \mapsto F^1_{(\w,u)}   (\g_s)\in G,$$
where $\G(t,s)=\g_s(t);s,t \in [0,1]$, is smooth. This follows from the fact {that} the map $t \in [0,1] \mapsto {\mathcal{H}}_\w(\g_s,t,u)$ is a solution of a differential equation in $t$. This differential equation depends smoothly on $s$ if $\G$ is smooth. For details see \cite[page 74]{KN}. 

The fact that $F^1_{(\w,u)}$ descends to a map $\pi_1^1(M,*)\to G$ follows from lemma \ref{EEE} and the fact that the horizontal lift of vectors defines a linear map. Finally the fact that $F^1_{(\w,u)}$ is a group morphism follows from the fact that {${\mathcal{H}}_\w\big (\g^2\g^1,1,u\big )={\mathcal{H}}_\w\big(\g^1,1,{\mathcal{H}}_\w(\g^2,1,u)\big)$,
for any smooth curves $\g^1$ and $\g^2$ in $M$ starting and ending at $*$.}\end{Proof}

It is proved in \cite{CP} that any holonomy $\pi_1^1(M,*) \to G$ arises from a bundle with connection in this way.

\subsubsection{Definition of a categorical holonomy}

\begin{Definition}[3-path]\label{K}
Let $M$ be a  manifold. 
A smooth map $J \colon[0,1]^3\to M$, say $J(t,s,x)=\G^x(t,s)=\g^x_s(t), $ where $t,s,x \in [0,1]$, is said to be a 3-path if:
\begin{enumerate}
\item $J$ is constant over $\{0\} \times [0,1]^2$ and over $\{1\} \times [0,1]^2$.
\item $J$ {restricts} to a rank-1 homotopy (connecting $\g^0_0$ with $\g^1_0$) over $[0,1] \times \{0\} \times [0,1]$ and over $[0,1] \times \{1\} \times [0,1]$ (connecting  $\g^0_1$ with $\g^1_1$).
\item $J$ has a product structure close to the boundary of $[0,1]^3$; see definition \ref{thin2}.
\end{enumerate}
\end{Definition}
In particular $J$ defines a homotopy connecting the 2-paths $\G^0$ and $\G^1$; see \ref{2es}. Note that it is not necessarily a rank-2 homotopy.
\begin{Definition}
Let $M$ and $N$ be smooth manifolds. Consider the projection map $S_2(M) \to \S_2(M) $. A map $F\colon \S_2(M) \to N$ is said to be smooth if, for any 3-path $J$,  the map $$x \in [0,1] \mapsto F([\G^x]) \in N$$  is smooth. Here $J(t,s,x)=\G^x(t,s)$, for each $(t,s,x) \in[0,1]^3 $.
\end{Definition}

Let \hbox{${\Gc= ( \d\colon E \to  G,\tr)}$} be a Lie crossed module. Consider the categorical group ${\mathcal C}(\Gc)$ constructed from $\Gc$, see \ref{CCG}. Its set of objects is $G$ and its set of morphisms is the semidirect product $G \ltimes E$ of  $G$ and $E$; all of these are smooth manifolds.
The following  definition extends the definition of holonomy for principal $G$-bundles of \cite{CP} and should be compared with analogous non-strict constructions in {\cite{BS,MP,M}.}

\begin{Definition}[{Categorical holonomy}]
Let $M$ be a smooth manifold. Choose an element $* \in M$.
A categorical holonomy is, {by definition,} a strict monoidal functor (in other words a categorical group morphism) $F\colon {\mathcal P}_{{2}}(M,*) \to \C(\Gc)$, such that the associated maps $F_1 \colon \pi_1^1(M) \to G$ and $F_2 \colon \P_2(M,*) \to G \ltimes   E$ on the sets of objects and morphisms of $\Pc_2(M,*)$ are smooth. Here ${\mathcal P}_{{2}}(M,*)$ is the thin fundamental categorical group of $M$; see \ref{TCM}.
\end{Definition}

If we use the notation 
$$[\G]=\begin{CD} \d_u([\G]) \\ @AA [\G]A \\   \d_d([\G]) \end{CD}\quad \quad \quad \textrm{   and } \quad \quad \quad  (X,e)=\begin{CD} \d(e)^{-1}X\\ @AA eA \\   X\end{CD}, $$
to denote the morphisms of $\Pc_2(M,*)$ and of $\C(\Gc)$, then the tensor product of two elements  $[\G]$ and $[\G']$ of $\P_2(M,*)$ can be presented in the  suggestive form:
\begin{equation} \begin{CD} \d_u([\G]) \\ @AA [\G]A \\   \d_d([\G]) \end{CD} \quad \tn \quad  \begin{CD} \d_u([\G']) \\ @AA [\G']A \\   \d_d([\G']) \end{CD} \quad = \quad \begin{CD} \d_u([\G]) \d_u([\G']) \\ @AA [\G] \circ_h [\G']A \\   \d_d([\G])  \d_d([\G']) \end{CD} 
\end{equation}
In this notation, the condition that a functor $F=(F_1,F_2)\colon \Pc_2(M,*) \to  \C(\Gc)$ be a strict monoidal functor means that:
\begin{equation}
 F\left ( \begin{CD} \d_u([\G]) \\ @AA [\G]A \\   \d_d([\G]) \end{CD} \quad \tn\quad \begin{CD} \d_u([\G']) \\ @AA [\G']A \\   \d_d([\G']) \end{CD} \right)=\begin{CD} F_1(\d_u([\G])) \\ @AA F_2([\G])A \\   F_1(\d_d([\G])) \end{CD} \quad \tn \begin{CD} F_1(\d_u([\G'])) \\ @AA F_2([\G'])A \\   F_1(\d_d([\G']) )\end{CD}
\end{equation}
for each $[\G],[\G'] \in \P_2(M,*) $. In other words:
\begin{equation}
 \begin{CD} F_1\big (\d_u([\G]) \d_u([\G'])\big) \\ @AA F_2\left([\G] \circ_h [\G'] \right) A \\   F_1\big (\d_d([\G])  \d_d([\G']) \big)\end{CD} \quad   \quad= \quad  \begin{CD} F_1\big (\d_u([\G])) F_1\left( \d_u([\G'])\right) \\ @AA \Big (F_1\left(\d_d[\G] \right)\tr  F_2\left( [\G'] \right)\Big ) F_2\left( [\G] \right) A \\   F_1\left (\d_d([\G])\right) F_1\left(  \d_d([\G']) \right)\end{CD} \quad  \quad  \quad  \quad  \quad 
\end{equation}
for each $[\G],[\G'] \in \P_2(M,*) ;$ see \ref{CCG}. Recall that an arrow $g \ra{(g,e)} \d(e^{-1}) g$ in $\C(\Gc)$ is sometimes simply denoted by $g \ra{e} \d(e^{-1}) g,$ where $ g \in G$ and $e \in E.$

Since the category of crossed modules and the category of  categorical groups are equivalent, we can give an alternative definition of a categorical holonomy, and define it as being a smooth crossed module map $$\Pi_2^{2,1}(M,*) \to \Gc,$$
where the notion of a smooth crossed module map is the obvious one.  Here $\Pi_2^{2,1}(M,*)$ is the thin fundamental crossed module of the smooth manifold $M$; see \ref{fcm}.

\subsection{Categorical connections and categorical connection holonomies}\label{ccch}
The aim of this subsection is to prove the following main theorem.

\begin{Theorem}\label{FFF}
Let $G$ be a Lie group. Let $M$ be a smooth manifold. Let also $P \to M$ be a smooth principal fibre bundle with structure group $G$. Let ${\Gc= ( \d\colon E \to  G,\tr)}$ be a Lie crossed module.
A $\Gc$-categorical connection $(\w,m)$ on $P$ together with an element $u\in P_*$, the fibre of  $P$ over $*\in M$, determines a categorical holonomy $F_{(\w,m,u)}\colon \Pc_2(M,*) \to \C(\Gc)$, which we call a categorical connection holonomy, since it arises from a categorical connection.
\end{Theorem}

{For the rest of this subsection, fix a Lie crossed module ${\Gc= ( \d\colon E \to  G,\tr)}$, and let ${\mathfrak{G}=(\d \colon \le \to  \lg,\tr )}$ be the associated differential crossed module. Consider also a  smooth manifold $M$ and a principal $G$-bundle $P$ over $M$, as well as a $\Gc$-categorical connection $(\w,m)$ on $P$. }

\subsubsection{Defining a categorical holonomy from a categorical connection}\label{dchcc}
Consider a base point $*\in M$. Choose also $u \in P_*$, the fibre of $P$ at $*$.
To describe $F_{(\w,m,u)}$ on the set $\pi_1^1(M,*)$ of objects of $\Pc_2(M,*)$, we use  {theorem} \ref{E}. Therefore we define:
\begin{equation}
F_{(\w,m,u)}^1([\g])=F_{(\w,u)}^{{1}}([\g]), \textrm{ where } [\g] \in \pi_1^1(M,*).
\end{equation}

Let now $\G$ be a 2-path in $M$ with $\d_l(\G)=\d_r(\G)=*$. Let $\g_s(t)=\G(t,s); s,t \in [0,1]$. To define $F^2_{(\w,m,u)}$, consider the smooth function $s \in [0,1] \mapsto e_\G(s)\in E$, which solves  the following differential equation in $E$:
\begin{equation}\label{G}
{\frac{d}{d s} e_\G(s)=e_\G(s) \int_{0}^{1} m \left (\t{\frac{\d}{\d t}\g_s(t)},\t{\frac{\d}{\d s}\g_s(t)} \right)  _{{\mathcal{H}}_\w(\g_s,t,u) } dt,}
\end{equation}
{with initial condition $e_\G(0)=1_E.$
Set
$e_\G=e_\G(1)$
to be the element of $E$ assigned to the 2-path $\G$.
Then we define:}
\begin{equation}
F^2_{(\w,m,u)}(\G)=\left ( F_{(\w,u)}^{{1}}\big ({[\d_d(\G)]}\big) ,e_\G\right).
\end{equation}

{Recall  the elements $g_{\g_s}\in G$, where $s \in [0,1]$, which are given by the condition:
$ug_{\g_s}={\mathcal{H}}_\w(\g_s,1,u); $
see \ref{BBB}.
Thus $F^1_{(\w,m,u)}({[\g_s]})=g^{-1}_{\g_s},$
where $g^{-1}_{\g_s}\doteq \big (g_{\g_s}\big )^{-1} $.  
 {In local coordinates, our definition of a categorical holonomy coming from a categorical connection coincides with definition 2.32 of \cite{BS}.}}

We will prove that  $F_{(\w,m,u)}=\left ( F^1_{(\w,m,u)},F^2_{(\w,m,u)}\right)$ defines a categorical holonomy  $\Pc_2(M,*) \to \C(\Gc)$. Explicitly, we need to prove that the assignment:

$$\begin{CD} \d_u(\G )\\ @AA\G A \\ \d_d(\G) \end{CD} \mapsto \begin{CD} g^{-1}_{\g_1}\\  @AA e_\G  A \\ g^{-1}_{\g_0} \end{CD} \quad \quad  ,$$
where $\G$ is a 2-path {in $M$} such that $\d_l(\G)=\d_r(\G)=*$, and $\G(t,s)=\g_s(t);\forall s,t \in [0,1]$, 
 depends only on the rank-2 homotopy class to which $\G$ belongs, and defines a monoidal functor $\Pc_2(M,*)\to \C(\Gc)$, such that the induced maps on the set of objects and morphisms of $\Pc_2(M,*)$ are smooth.

First of all:
\begin{Lemma}\label{DDD}
{Let  $\G$ be a 2-path in $M$ with $\d_l(\G)=\d_r(\G)=*$. Then  for each $s \in [0,1]$ we have:
$ g_{\g_0}\d(e_\G(s))= g_{\g_s}.$
Here $\G(t,s)=\g_s(t);\forall s,t \in [0,1]$.}
\end{Lemma}
\begin{Proof}
{The previous equation holds for $s=0$, by definition of $e_\G(s)$. We have:}
$$
\frac{d}{d s} \big( g_{\g_0}\d(e_\G(s)) \big)= g_{\g_0}\frac{d}{d s}\d\left ( e_\G(s)\right)= g_{\g_0}\d\left(\frac{d}{d s} e_\G(s)\right),$$
therefore:
 \begin{align*}\frac{d}{d s} \big( g_{\g_0}\d(e_\G(s)) \big)&= {g_{\g_0}\d  \Big( e_\G(s) \int_{0}^{1} m \left (\t{\frac{\d}{\d t}\g_s(t)},\t{\frac{\d}{\d s}\g_s(t)} \right)  _{{\mathcal{H}}_\w(\g_s,t,u) } dt \Big)} \\
&= g_{\g_0}\d  ( e_\G(s)) \int_{0}^{1} \W \left (\t{\frac{\d}{\d t}\g_s(t)},\t{\frac{\d}{\d s}\g_s(t)} \right)  _{{\mathcal{H}}_\w(\g_s,t,u) }  dt,
\end{align*}
since $\d(m)=\W$.
By equation (\ref{C}) we also have:
\begin{align*}
 \frac{d}{d s} g_{\g_s}    =       g_{\g_s}              \int_{0}^{1} \W \left (\t{\frac{\d}{\d t}\g_s(t)},\t{\frac{\d}{\d s}\g_s(t)} \right)  _{{\mathcal{H}}_\w(\g_s,t,u) }.
    \end{align*}
Thus both sides of the equation are solutions of the same differential equation in $G$, and have  the same initial condition.
\end{Proof}

{We remark that, using the terminology of \cite{BS}, this result corresponds to the condition that a local pre-$2$-holonomy gives rise to a local true $2$-holonomy, which there also follows from the vanishing of the fake curvature.}

Let $\G \colon [0,1]^2 \to M$ be a 2-path in $M$ with $\d_l(\G)=\d_r(\G)=*$. 
We   need to prove that $F^2_ {(\w,m,u)} (\G)$ depends only on the 2-track defined by $\G$, in other words on the rank-2 homotopy class of 2-paths to which $\G$ belongs. This will be shown shortly in Corollary \ref{cormain2} of Theorem \ref{Main2}. For now, 
 consider a 3-path $J$, where $J(t,s,x)=\G^x(t,s); \forall t,s,x \in [0,1]$, defining {a homotopy} (not necessarily rank-2) connecting the 2-paths $\G^0$ and $\G^1$. Suppose that $\d_l(\G^0),\d_r(\G^0)=*$, thus $\d_l(\G^1),\d_r(\G^1)=*$. We want to calculate {$\frac{\d}{\d x}e_{\G^x}(s),$}
where {$s,x \in [0,1]$.}

In fact, for later applications, we will consider the more general case when, for each $x \in [0,1]$, we have that {$\d_l(\G^x)=q(x)$ and $\d_r(\G^x)=q(x)$,} where $q\colon [0,1] \to M$ is a smooth map. {This generality will be important to define Wilson spheres in \ref{WS}.} The lemma in the following subsection will be useful for our purpose.

\subsubsection{A well-known lemma}\label{wkl}
Let $G$ be a Lie group. Consider a $\lg$-valued smooth function $V(s,x)$ defined on $[0,1]^2$.
Consider  the following differential equation in $G$:
 $$\frac{\d}{\d s} a(s,x)=a(s,x)V(s,x),$$ 
with initial condition $a(0,x)=1_G, \forall x \in [0,1]$.
We want to know $\frac{\d}{\d x} a(s,x)$.

Let $\theta$ be the canonical $\lg$-valued $1$-form on $G$. Thus $\theta$ is left {invariant} and satisfies $\theta(A)=A, \forall A \in \lg$, being defined uniquely by these properties. Also 
\begin{equation}\label{AAA}
d \theta(A,B)=- \theta ([A,B]),
\end{equation}
where $A,B \in \lg$. 
We have:
$$\frac{\d }{ \d x} \theta\left (\frac{\d}{\d s} a(s,x) \right )= \frac{\d }{ \d x}\theta\big (  a(s,x)V(s,x) \big )= \frac{\d }{ \d x}V(s,x).
$$
On the other hand: 
\begin{align*}\frac{\d }{ \d x} \theta\left (\frac{\d}{\d s} a(s,x) \right ) 
&= d a^*(\theta)\left(  \frac{\d }{ \d x}, \frac{\d }{ \d s}\right) {+}\frac{\d }{ \d s} a^*(\theta) \left ( \frac{\d }{ \d x} \right){+} a^*(\theta)\left(  \left [\frac{\d }{ \d x}, \frac{\d }{ \d s} \right]\right) \\
&=d \theta\left(  \frac{\d }{ \d x}  a(s,x) , \frac{\d }{ \d s} a(s,x) \right)  {+}\frac{\d }{ \d s} \theta \left ( \frac{\d }{ \d x}a(s,x) \right).
\end{align*}  
Therefore:
\begin{multline*}\theta \left ( \frac{\d }{ \d x} a(s,x) \right)\\= \int_0^s \left ({-}d \theta\left(  \frac{\d }{ \d x}  a(s',x) , \frac{\d }{ \d s'} a(s',x) \right) {+}\frac{\d }{ \d x} V(s',x)\right )ds'+\theta\left (\frac{\d }{ \d x}a(0,x)\right).
\end{multline*}
Since  $\frac{\d }{ \d x}a(0,x)=0$ (due to the initial conditions) we have the following:

\begin{Lemma}\label{WKL}
$$ \frac{\d }{ \d x} a(s,x)= a(s,x)  \int_0^s \left ({-}d \theta\left(  \frac{\d }{ \d x}  a(s',x) , \frac{\d }{ \d s'} a(s',x) \right) {+} \frac{\d }{ \d x} V(s',x)\right )ds',$$
{for each $x,s \in [0,1]$.}
\end{Lemma}

\subsubsection{The dependence of the surface holonomy on a smooth family of surfaces}

The discussion in this subsection will be very similar to that of \ref{BBB}. We want to prove an analogue of lemma \ref{EEE} for the surface holonomy of a family of 2-paths {$\G^x$.}

Consider a smooth map $J\colon [0,1]^3 \to M$, where   $J(t,s,x)=\G^x(t,s); \forall t,s,x \in [0,1]$.
 Given $x,s \in [0,1]$, a smooth curve  $\g^x_s\colon [0,1] \to M$ is defined as $\g^x_s(t)=J(t,s,x)$, where $t,s,x \in [0,1]$.  Suppose that, for each $x \in [0,1]$, there exists $q(x) \in M$ such that $J(0,s,x)=J(1,s,x)=q(x)$ for each $s \in  [0,1]$. In other words $\d_l(\G^x)=\d_r(\G^x)=q(x)$.
Note that the map $x \in [0,1] \mapsto q(x) \in M$ is necessarily smooth. Let $q(0)=* \in M$. Choose $ u \in P_*$, the fibre of $P$ at $*$. Let also $u_x={\mathcal{H}}_\w(q,x,u)$, where $x \in [0,1]$. In particular $u_0=u$.  

Note that the conditions above imply that:
\begin{equation}\label{S1}
\Rank ({ {\EuScript{D}}}_{(s,x)} J)_{(t,s,x)}\leq 1 \textrm{ if } t=0 \textrm{ or } t=1,\textrm{ for each } s,x \in [0,1].
\end{equation}
Suppose that we also  have:
\begin{equation}\label{S2}
\Rank ({{\EuScript{D}}}_{(t,x)} J)_{(t,s,x)}\leq 1 \textrm{ if } s=0 \textrm{ or } s=1,\textrm{ for each } t,x \in [0,1]. 
\end{equation}

For each $x \in [0,1]$, let $s \in [0,1] \mapsto e_{\G^x}(s)\in E$ be the solution of the differential equation in $E$:
\begin{equation}\label{E1}
{\frac{d}{d s} e_{\G^x}(s)=e_{\G^x}(s)\int_{0}^{1} m \left (\t{\frac{\d}{\d t}\g^x_s(t)},\t{\frac{\d}{\d s}\g^x_s(t)} \right)  _{{\mathcal{H}}_\w(\g_s^x,t,u_x) } dt ,}
\end{equation}
with initial condition:
\begin{equation}\label{E2}
e_{\G^x}(0)=1_E, \forall x \in [0,1].
\end{equation}
See equation (\ref{G}). Our purpose is to calculate the $x$-dependence of $e_{\G^x}= e_{\G^x}(1)$, the element of $E$ assigned to $\G_x$, by calculating
{${\frac{d}{d x}} e_{\G^x} .$}
The result of this calculation, Theorem \ref{Main2}, is entirely analogous to Lemma \ref{EEE}, with the 2-curvature 3-form $\M$ replacing the curvature 2-form $\W$.

Lemma (\ref{WKL})    leads to:
\begin{multline} \frac{d }{ d x} e_{\G^x} =e_{\G^x}  \int_0^1 \int_{0}^{1} \frac{\d }{ \d x}  \left (m \left (\t{\frac{\d}{\d t}\g^x_s(t)},\t{\frac{\d}{\d s}\g^x_s(t)} \right)  _{{\mathcal{H}}_\w(\g_{s} ^x,t,u_x) }  \right) dt ds\\
{-} e_{\G^x}  \int_0^1   d \theta\left(  \frac{\d }{ \d x}  e_{\G^x}(s) , \frac{\d }{ \d s} e_{\G^x}(s) \right)   ds.
\end{multline}
 Let
\begin{align*}
A_x&=\int_0^1 \int_{0}^{1} \frac{\d }{ \d x} \left ( m \left (\t{\frac{\d}{\d t}\g^x_s(t)},\t{\frac{\d}{\d s}\g^x_s(t)} \right)  _{{\mathcal{H}}_\w(\g_{s} ^x,t,u_x) } \right) dt ds\in \le\\
B_x&=\int_0^1 d \theta\left(  \frac{\d }{ \d x}  e_{\G^x}(s) , \frac{\d }{ \d s} e_{\G^x}(s) \right)   ds \in \le.
\end{align*}
Thus 
\begin{equation}\label{torefer}
 \frac{d }{ d x} e_{\G^x} =e_{\G^x} (A_x-B_x).
\end{equation}

Let us analyse $A_x$ and $B_x$ separately.
Consider the map $f\colon [0,1]^3 \to P$ such {that} $f(x,s,t)={{\mathcal{H}}_\w(\g_s^x,t,u_x) }$, for each $ x,s,t \in [0,1].$ This map is smooth; see \cite[page 74]{KN}. {By definition  we have 
$\frac{\d}{\d t } f(x,s,t) = \t{ \frac{\d}{\d t}\g^x_s(t)}_{{\mathcal{H}}_\w(\g_s^x,t,u_x) }.$
Also, we trivially get: $
\left (\frac{\d}{\d s } f(x,s,t)\right) ^H = \t{ \frac{\d}{\d s}\g^x_s(t)}_{{\mathcal{H}}_\w(\g_{s}^x,t,u_x) }$ and $
\left ( \frac{\d}{\d x } f(x,s,t)\right)^H = \t{ \frac{\d}{\d x}\g^x_s(t)}_{{\mathcal{H}}_\w(\g_{s}^x,t,u_x) }$;
see \ref{BBB}.}
Therefore, since $m(X,Y)$ vanishes if either $X$ or $Y$ is vertical, it follows that:
\begin{align*} 
A_x&= \int_{0}^{1}\int_{0}^{1} \frac{\d }{ \d x}  m \left (\frac{\d}{\d t} f(x,s,t),\frac{\d}{\d s}f(x,s,t) \right) dt ds\\ 
  &  =    \int_0^1 \int_{0}^{1} dm\left ( \frac{\d}{\d x} f(x,s,t), \frac{\d}{\d t}f(x,s,t), \frac{\d}{\d s}f(x,s,t)\right)dtds
\\ & \quad \quad\quad \quad\quad \quad {-} \int_0^1 \int_{0}^{1} \frac{\d}{\d t} m\left ( \frac{\d}{\d s}f(x,s,t), \frac{\d}{\d x}f(x,s,t)\right )dtds
\\  &\quad \quad\quad \quad \quad \quad {-}  \int_0^1 \int_{0}^{1} \frac{\d}{\d s} m\left ( \frac{\d}{\d x}f(x,s,t), \frac{\d}{\d t}f(x,s,t)\right )dtds.
\end{align*}
{We have used the well known equation:}
\begin{multline*}
d\a(X,Y,Z)=X\a(Y,Z)+Y\a(Z,X)+Z\a(X,Y)\\+\a(X,[Y,Z])+ \a(Y,[Z,X])+\a(Z,[X,Y]),
\end{multline*}
valid for any smooth 2-form {$\a$} in a manifold.
Therefore:
\begin{align*}
A_x&=\int_{0}^{1}\int_{0}^{1} dm\left ( \frac{\d}{\d x} f(x,s,t), \frac{\d}{\d t}f(x,s,t), \frac{\d}{\d s}f(x,s,t)\right)dtds
\\ &{-} \int_0^1  \left [m\left ( \frac{\d}{\d s}f(x,s,1), \frac{\d}{\d x}f(x,s,1)\right )-m\left ( \frac{\d}{\d s}f(x,s,0), \frac{\d}{\d x}f(x,s,0)\right )\right]ds
\\  &{-} \int_0^1  \left [m\left ( \frac{\d}{\d x}f(x,1,t), \frac{\d}{\d t}f(x,1,t)\right ) - m\left ( \frac{\d}{\d x}f(x,0,t), \frac{\d}{\d t}f(x,0,t)\right )\right]dt.
\end{align*}
{We analyse each term separately.} Since $m(X,Y)=0$ if either $X$ or $Y$ is vertical we have:
$${
m\left ( \frac{\d}{\d x}f(x,1,t), \frac{\d}{\d t}f(x,1,t)\right )=m\left (\t{\frac{\d}{\d x}\g^x_1(t)},\t{\frac{\d}{\d t}\g^x_1(t)} \right)  _{{\mathcal{H}}_\w(\g_1^x,t,u_x) } =0,}$$
by equation (\ref{S2}).
Analogously from equation  (\ref{S1}) or equation  (\ref{S2}) it follows that:
$m\left ( \frac{\d}{\d x}f(x,0,t), \frac{\d}{\d t}f(x,0,t)\right )=0,$
 $m\left ( \frac{\d}{\d s}f(x,s,0), \frac{\d}{\d x}f(x,s,0)\right )=0,$ and
$m\left ( \frac{\d}{\d s}f(x,s,1), \frac{\d}{\d x}f(x,s,1)\right )=0.$
Therefore we have:
\begin{equation}\label{torefer2}
A_x=\int_{0}^{1}\int_{0}^{1} dm\left ( \frac{\d}{\d x} f(x,s,t), \frac{\d}{\d t}f(x,s,t), \frac{\d}{\d s}f(x,s,t)\right)dtds,
\end{equation}
for each $x \in [0,1]$.

{We now analyse $B_x$, for each $x \in [0,1]$. We have:}
\begin{align*}
& d \theta\Big(  \frac{\d }{ \d x}  e_{\G^x}(s) , \frac{\d }{ \d s} e_{\G^x}(s) \Big)  
\\ &= \int_0^1   d \theta\Big(  \frac{\d }{ \d x}  e_{\G^x}(s) ,e_{\G^x}(s) m\left (\t{\frac{\d}{\d t}\g^x_s(t)},\t{\frac{\d}{\d s}\g^x_s(t)} \right)_{{\mathcal{H}}_\w(\g_s^x,t,u_x) }  \Big)  dt,\\
&= \int_0^1   d \theta\Big( e_{\G^x}(s)^{-1} \frac{\d }{ \d x}  e_{\G^x}(s) , m \left (\t{\frac{\d}{\d t}\g^x_s(t)},\t{\frac{\d}{\d s}\g^x_s(t)} \right )_{{\mathcal{H}}_\w(\g_s^x,t,u_x) } \Big )dt
\\& = -\int_0^1   \Big[ e_{\G^x}(s)^{-1} \frac{\d }{ \d x}  e_{\G^x}(s) , m \left (\t{\frac{\d}{\d t}\g^x_s(t)},\t{\frac{\d}{\d s}\g^x_s(t)}  \right)	 _{{\mathcal{H}}_\w(\g_s^x,t,u_x) }\Big]   dt
\\& = -\int_0^1  \partial \left ( e_{\G^x}(s)^{-1} \frac{\d }{ \d x}  e_{\G^x}(s) \right ) \triangleright m \left (\t{\frac{\d}{\d t}\g^x_s(t)},\t{\frac{\d}{\d s}\g^x_s(t)} \right)	 _{{\mathcal{H}}_\w(\g_s^x,t,u_x) }   dt.
\end{align*}
The first equation follows by definition of $ e_{\G^x}(s)$, the second since the $\le$-valued 2-form $d \theta$ on $E$ is left {invariant} (since $\theta$ is left {invariant}) and the third by the well-known formula (\ref{AAA}). Finally, the last equation  follows since $(\partial\colon \le \to \lg, \tr )$ is a differential crossed module.

{By lemma \ref{DDD}, we also  have that  $g_{\g^x_0}\d(e_{\G^x}(s))= g_{\g^x_s}$, for each $s,x \in [0,1]$. Therefore:}
\begin{align*}
\partial \Big ( e_{\G^x}^{-1}(s) \frac{\d }{ \d x}  e_{\G^x}(s) \Big )&=
   \partial \Big ( e_{\G^x}^{-1} (s)\Big) \partial \Big(\frac{\d }{ \d x}  e_{\G^x}(s) \Big )\\
&= g_{\g^x_s}^{-1} g_{\g^x_0}     \frac{\d }{ \d x} \Big(g_{\g^x_0}^{-1} g_{\g^x_s}\Big)         \\
&= \int_{0}^{1} \W \left (\t{\frac{\d}{\d t}\g_s^x(t)},\t{\frac{\d}{\d x}\g_s^x(t)} \right) _{{\mathcal{H}}_\w(\g_s^x,t,u_x) } dt,
\end{align*}
where the last equation follows  by lemma \ref{EEE}.  We also use the fact that  $\frac{\d }{ \d x} g_{\g^x_0}=0.$ This follows from lemma \ref{EEE} and equation \ref{S2}. Putting everything together, we obtain:
\begin{multline*}
{d \theta\left(  \frac{\d }{ \d x}  e_{\G^x}(s) , \frac{\d }{ \d s} e_{\G^x}(s) \right)}
 =-\int_{0}^{1} \W \left (\t{\frac{\d}{\d t}\g_s^x(t)},\t{\frac{\d}{\d x}\g_s^x(t)} \right) _{{\mathcal{H}}_\w(\g_s^x,t,u_x) } dt\,   \tr\\ \int_0^1 m \left (\t{\frac{\d}{\d t}\g^x_s(t)},\t{\frac{\d}{\d s}\g^x_s(t)}  \right) _{{\mathcal{H}}_\w(\g_s^x,t,u_x) }   dt.
\end{multline*}
{Due to the fact that  $(\d\colon \le  \to \lg, \tr )$  is a differential crossed module,  given $e,f \in \le$ we have:
$\d(e)\tr f=[e,f]=-[f,e]=-\d(f) \tr e.$
If we use this  together with the fact that  $\W(X,Y)$ and $m(X,Y)$  vanish if either of the vectors $X$ or $Y$ is vertical it thus follows that (since $\d(m)=\W$):}
\begin{multline}\label{tcompare}
{B_x} =
\int_0^1 \left (\int_{0}^{1} \W \left (\frac{\d}{\d t}f{(x,s,t)},\frac{\d}{\d s}f{(x,s,t)}  dt \right )  \tr\right. \\ \left .\int_0^1 m \left (\frac{\d}{\d t}f{(x,s,t)},\frac{\d}{\d x}f{(x,s,t)}  \right)   dt\right)ds.
\end{multline}

The following is a 2-dimensional version of Lemma \ref{EEE}.
\begin{Theorem}\label{Main2}
Let $M$ be a smooth manifold. Let ${\Gc=(\d\colon E \to G,\tr)}$ be a Lie crossed module. Let $P \to M$ be a principal $G$-bundle over $M$. Consider a $\Gc$-categorical connection $(\w,m)$ on $P$.
Let $J\colon [0,1]^3 \to M$ be a smooth map. Let $J(t,s,x)=\G^x(t,s)=\g^x_s(t);\forall t,s,x \in [0,1]$. Suppose that $\d_l(\G^x)=\d_r(\G^x)=q(x),\forall x \in [0,1]$, where $q\colon [0,1] \to M$ is a smooth curve. Choose $u \in P_{q(0)}$, the fibre of $P$ at $q(0)$. Let also $u_x=\H_\w(q,{x},u)$, where $x \in [0,1]$. In particular $u_0=u$. 

Consider the map $(s,x)\in [0,1]^2 \mapsto e_{\G^x}(s)\in E$ defined by equations (\ref{E1}) and  (\ref{E2}). {Let $e_{\G^x}=e_{\G^x}(1)$}. For each $x\in [0,1]$, we have:

\begin{align*}\frac{d }{ d x}e_{\G^x}  &=e_{\G^x} \int_0^1 \int_{0}^{1} dm\left ( \t{\frac{\d}{\d x} \g_s^x(t)},\t{ \frac{\d}{\d t} \g_s^x(t)},\t{ \frac{\d}{\d s} \g_s^x(t)}\right)_{{\mathcal{H}}_\w(\g_s^x,t,u_x) }  dtds\\
&=e_{\G^x} \int_0^1 \int_{0}^{1} \M\left ( \t{\frac{\d}{\d x} \g_s^x(t)},\t{ \frac{\d}{\d t} \g_s^x(t)},\t{ \frac{\d}{\d s} \g_s^x(t)}\right)_{{\mathcal{H}}_\w(\g_s^x,t,u_x) }  dtds.
\end{align*}
\end{Theorem}
\begin{Proof}
By the 2-structure equation, see \ref{se}{, and equation (\ref{BB})} it follows that:
\begin{align*}
\int_0^1 \int_{0}^{1} &dm\left ( \t{\frac{\d}{\d x} \g_s^x(t)},\t{ \frac{\d}{\d t} \g_s^x(t)},\t{ \frac{\d}{\d s} \g_s^x(t)}\right)_{{\mathcal{H}}_\w(\g_s^x,t,u) }  dtds\\
&=\int_0^1 \int_{0}^{1} \M\left ( \frac{\d}{\d x} f{(x,s,t)}, \frac{\d}{\d t}f{(x,s,t)}, \frac{\d}{\d s} f{(x,s,t)}\right) dtds\\
&=\int_0^1 \int_{0}^{1} dm\left ( \frac{\d}{\d x} f{(x,s,t)}, \frac{\d}{\d t}f{(x,s,t)}, \frac{\d}{\d s} f{(x,s,t)}\right) dtds\\
& \quad \quad +\int_0^1 \int_{0}^{1} \w\left ( \frac{\d}{\d x} f{(x,s,t)}\right) \tr m\left (\frac{\d}{\d t}f{(x,s,t)}, \frac{\d}{\d s} f{(x,s,t)}\right) dtds\\
& \quad \quad -\int_0^1 \int_{0}^{1} \w\left ( \frac{\d}{\d s} f{(x,s,t)}\right) \tr m\left (\frac{\d}{\d t}f{(x,s,t)}, \frac{\d}{\d x} f{(x,s,t)}\right) dtds.
\end{align*}
Note     that   $\w\left ( \frac{\d}{\d t} f(x,t,s\right )=\w\left (\frac{\d}{\d t}{\mathcal{H}}_\w(\g^x_s,t,u_x)\right) =0$.

{Also by {lemma} \ref{Main1} and the fact that $\W$ is horizontal:}
\begin{align*}
&\int_{0}^{1} \w\left ( \frac{\d}{\d s} f{(x,s,t)}\right) \tr m\left (\frac{\d}{\d t}f{(x,s,t)}, \frac{\d}{\d x} f{(x,s,t)}\right) dt \\
   &= \int_{0}^{1} \int_0^{t} \W  \left (\frac{\d}{\d t'}f{(x,s,t')}, \frac{\d}{\d s} f{(x,s,t')}\right){dt'} \tr m\left (\frac{\d}{\d t}f{(x,s,t)}, \frac{\d}{\d x} f{(x,s,t)}\right) dt \\
&= \int_{0}^{1}   \W  \left (\frac{\d}{\d t'}f{(x,s,t')}, \frac{\d}{\d s} f{(x,s,t')}\right)dt' \tr \int_0^1 m\left (\frac{\d}{\d t'}f{(x,s,t')}, \frac{\d}{\d x} f{(x,s,t')}\right) dt' \\          
&  - \int_0^1\W\left (\frac{\d}{\d t}f{(x,s,t)}, \frac{\d}{\d s} f{(x,s,t)}\right) \tr \left (\int_0^t m\left (\frac{\d}{\d t'}f{(x,s,t')}, \frac{\d}{\d x} f{(x,s,t')}\right) dt' \right )dt.
\end{align*}
The last equation follows from integrating by parts. Indeed, given a smooth $\lg$-valued function $V(t)$ and a smooth $\le$-valued function $W(t)$, each defined on $[0,1]$, we have:
\begin{multline*} \int_0^1 \int_0^tV(t')dt' \tr W(t)dt\\=\left [ \left (\int_0^tV(t')dt'\right) \tr \left (\int_0^tW(t')dt'\right)\right]_{t=1}- \int_0^1 V(t) \tr \left(\int_0^tW(t')dt'\right) dt.\end{multline*}
{Let us analyse the very last term. Applying the fact that  $(\d\colon \le \to \lg, \tr  )$ is a differential crossed module together with $\d(m)=\W$, part of the definition of a categorical connection $(\w,m)$, we have, since $\d(e)\tr f=-\d(f) \tr e; \forall e,f \in {\le}$:}
\begin{align*}
&\int_0^1\W\left (\frac{\d}{\d t}f{(x,s,t)}, \frac{\d}{\d s} f{(x,s,t)}\right) \tr \left (\int_0^t m\left (\frac{\d}{\d t'}f{(x,s,t')}, \frac{\d}{\d x} f{(x,s,t')}\right) dt' \right )dt\\
 &=-\int_0^1 \int_0^t \W\left (\frac{\d}{\d t'}f{(x,s,t')}, \frac{\d}{\d x} f{(x,s,t')}\right) dt' \tr m\left (\frac{\d}{\d t}f{(x,s,t)}, \frac{\d}{\d s} f{(x,s,t)}\right)   dt\\
&=- \int_0^1 \w\left (\frac{\d}{\d x}f{(x,s,t)}\right) dt' \tr m\left (\frac{\d}{\d t}f{(x,s,t)}, \frac{\d}{\d s} f{(x,s,t)}\right)   dt.
\end{align*}
{Therefore:}
$$
\int_0^1 \int_{0}^{1} \M\left ( \t{\frac{\d}{\d x} \g_s^x(t)},\t{ \frac{\d}{\d t} \g_s^x(t)},\t{ \frac{\d}{\d s} \g_s^x(t)}\right)_{{\mathcal{H}}_\w(\g_s^x,t,{u_x}) }  dtds =A_x-B_x,
$$
by equations (\ref{torefer2}) and (\ref{tcompare}). Comparing with equation (\ref{torefer}), yields theorem \ref{Main2}.
\end{Proof}

From this theorem and the fact that the horizontal lift $X \mapsto \t{X}$ of vector fields on $M$ defines a linear map $\X(M) \to \X(P)$ we obtain the following:
\begin{Corollary}\label{cormain2}
Let $M$ be a smooth manifold with a base point $*$. Let also ${\Gc=(\d\colon E \to G,\tr)}$ be a Lie crossed module. Let $P \to M$ be a principal $G$-bundle over $M$, and let $u \in P_*$, the fibre at $*\in M$. Consider a $\Gc$-categorical connection $(\w,m)$ on $P$. Let $\G$ be a 2-path in $M$ with $\d_l(\G)=\d_r(\G)=*$. Let $\G(t,s)=\g_s(t)$, for each $ t,s \in [0,1]$. Then $F_{(\w,m,u)}^2(\G)=(g_{\g_0}^{-1},e_\G )$ depends only on the rank-2 homotopy class of 2-tracks $[\G]$ to which $\G$ belongs.
\end{Corollary}

Consider a  crossed module ${\Gc=(\d\colon E \to G,\tr)}$. Let the categorical group it defines be $\C(\Gc)$; see \ref{CCG}. Let also be given a smooth manifold $M$ with a base point $*$, a principal $G$-bundle $P$ over $M$ and $\Gc$-categorical connection $(\w,m)$ on $P$.  We  defined maps $F^1_{(\w,m,u)}\colon \pi_1^1(M,*) \to G$  (the sets of objects of $\Pc(M,*)$ and of  $\C(\Gc)$, respectively) and $F^2_{(\w,m,u)}\colon \P_2(M,*) \to G \ltimes E$ (the sets of {morphisms} of $\Pc(M,*)$ and   $\C(\Gc)$, respectively). We now prove that this assignment, defined in \ref{dchcc}, is a categorical holonomy. This will complete the proof of theorem \ref{FFF}.

\subsubsection{Re-scaling}
Let $\G$ be a 2-path in $M$ such that $\d_l(\G)=\d_r(\G)=*$. As usual put $\G(t,s)=\g_s(t); s,t \in [0,1]$. Consider a smooth map $f \colon [a,b] \to [0,1]$, where $-\infty < a <b <+\infty$. We suppose that $f(a)=0$ and $f(b)=1$. Let $\G'(s,t)=\g'_s(t)$ be such that  $\g'_s(t)=\g_{f(s)}(t)$, where $s \in [a,b]$ and $t \in [0,1]$. Even though $\G'$ is not a 2-path, for its domain is $[0,1]\times [a,b]$, the elements $g_{\g'_s}\in E$ and $e_{\G'}(s)\in E$, where $s \in [a,b]$, defined in \ref{dchcc}, still make sense.
\begin{Lemma}\label{rescaling}
We have 
$${g_{\g'_s}=g_{\g_{f(s)}} \quad \textrm{and}\quad
e_{\G'}(s)=e_{\G}(f(s)),}
$$
for each $s\in [a,b]$.
\end{Lemma}
\begin{Proof}
The first equation follows since $\g'_s=\g_{f(s)}$, as curves $[0,1] \to M$, for each $s\in [a,b]$. As for the second one we have:
\begin{align*}
\frac{d}{d s} e_{\G'}(s)&={e_{\G'}(s) \int_{0}^{1} m \left (\t{\frac{\d}{\d t}{\G'}(t,s)},\t{\frac{\d}{\d s}{\G'}(t,s)} \right)_{{\mathcal{H}}_\w(\g'_s,t,u) } dt}\\
&{=e_{\G'}(s)\frac{d}{d s} f(s)\int_{0}^{1} m \left (\t{ \frac{\d \G }{\d t}(t,f(s))},\t{\frac{\d \G}{\d s}(t,f(s))} \right)  _{{\mathcal{H}}_\w(\g_{f(s)},t,u) } dt .}
\end{align*}
Whereas:
$$
{\frac{d}{d s} e_\G(f(s))= e_\G(f(s))\frac{d}{d s} f(s)\int_{0}^{1} m \left (\t{ \frac{\d \G }{\d t}(t,f(s))},\t{\frac{\d \G}{\d s}(t,f(s))} \right)  _{{\mathcal{H}}_\w(\g_{f(s)},t,u) } dt .}
$$
Therefore $s \mapsto e_{\G'}(s)$ and $s \mapsto e_\G(f(s))$  are both solutions of the same differential equation and they have the same initial conditions.
\end{Proof}

\subsubsection{Verification of the axioms for a categorical holonomy}
Let ${\Gc=(\d\colon E \to G,\tr)}$ be a Lie crossed module and let ${\mathfrak{G}=(\d \colon \le \to  \lg,\tr )}$ be the associated differential crossed module. 
Fix a smooth manifold $M$ with a base point $*$, a principal $G$-bundle $P$ over $M$ and an element $u\in P_*$, the fibre of $P$ over $*$. Consider a categorical connection $(\w,m)$ on $P$. Let us prove that the assignments 
$$[\g] \in \pi_1^1(M,*) \mapsto F^1_{(\w,u)}([\g]) \in G$$
and
$$[\G] \in \P_2(M,*) \mapsto  F^2_{(\w,m,u)}([\G])  \in G \ltimes E,$$
defined in \ref{dchcc}, yield a categorical holonomy $F_{(\w,m,u)}\colon \Pc(M,*) \to \C(\Gc)$. To begin with, let us see that $F_{(\w,m,u)}$  defines a monoidal functor $\Pc_2(M,*) \to \C(\Gc)$. Here $\C(\Gc)$ is the  categorical group constructed from $\Gc$; see \ref{CCG}.

{\bf Claim:}
$F_{(\w,m,u)}$ is a functor.

\begin{Proof}
Recall that the source and target maps in  a category are denoted by $\sigma$ and $\tau$, respectively. Let $\G$ be  a 2-path in $M$ such that $\d_l(\G),\d_r(\G)=*$.
 As usual put $\G(t,s)=\g_s(t); t,s \in [0,1].$ We have:
$$
{\sigma\left (F^2_{(\w,m,u)}({[\G]})\right)=\sigma\left(g_{\g_0}^{-1},e_{{[\G]}} \right)= g_{\g_0}^{-1}={F_1\left(\sigma({[\G]})\right).}}$$
On the other hand by lemma \ref{DDD}:
$$
\tau\left (F^2_{(\w,m,u)}({[\G]})\right)=\tau\left(g_{\g_0}^{-1},e_{{[\G]}} \right)= \d(e_{{\G}}^{-1} )g_{\g_0}^{-1}= g_{\g_1}^{-1} ={F_1\left(\tau({[\G]})\right)}.$$
 Therefore we can represent  $F_{(\w,m,u)}$ as:
$$\begin{CD} \d_u([\G]) \\ @AA [\G]A \\   \d_d([\G]) \end{CD} \in \Pc_2(M,*) \quad \longmapsto \quad \begin{CD}  g_{\g_1}^{-1}= \d(e_\G )^{-1}g_{\g_0}^{-1}\\ @AA {e_\G} A \\   g_{\g_0}^{-1} \end{CD} \in \C(\Gc) , $$
where, as usual, $\G(t,s)=\g_s(t); s,t \in  [0,1]$.

Let $\G^1$ and $\G^2$  be   2-paths in $M$ with $\d_l(\G^1),\d_r(\G^1)=*$ and $\d_l(\G^2),\d_r(\G^2)=*$. Suppose also that $\d_u(\G^1)=\d_d(\G^2)$. For each $s,t \in [0,1]$, put $\G^1(t,s)=\g_s^1(t)$  and $\G^2(t,s)=\g_s^2(t)$. Let $\G=\G_1\circ_v \G_2$. Put 
$\G(t,s)=\g_s(t); s,t \in [0,1]$.
Suppose $s \in [0,1/2]$. We have:
$$
{\frac{d}{d s} e_\G(s)=e_\G(s)\int_{0}^{1} m \left (\t{\frac{\d}{\d t}\g_s(t)},\t{\frac{\d}{\d s}\g_s(t)} \right)  _{{\mathcal{H}}_\w(\g_s,t,u) } dt.}
$$
By lemma \ref{rescaling}, $e_\G(1/2)=e_{\G^1} $. Let $e'_\G(s)= e_\G(1/2)^{-1} e_\G(s)$, where $s \in [1/2,1]$.  Then:
$$
{\frac{d}{d s} e'_\G(s) =e'_\G(s)\int_{0}^{1} m \left (\t{\frac{\d}{\d t}\g_s(t)},\t{\frac{\d}{\d s}\g_s(t)} \right)  _{{\mathcal{H}}_\w(\g_s,t,u) } dt.}
$$
Since {$e'_\G(1/2)=1_E$,} it follows, by lemma \ref{rescaling}, that $e'_\G =e_{\G^2} $.  Putting everything together yields $e_{\G} =e_{\G^1} e_{\G^2} $. Therefore we have:
$$F^2_{(\w,m,u)}\left ( \begin{CD} \d_u([\G^2] )\\ @AA[\G^2] A\\  \d_d([\G^2])\\ @AA[\G^1] A\\ \d_d([\G^1]) \end{CD}\right)\quad =\quad  \begin{CD} F^1_{(\w,m,u)}\left (\d_u([\G^2] )\right)\\ @AAF^2_{(\w,m,u)}\big ([\G^2]\big) A\\ F^1_{(\w,m,u)}\big( {\d_u([\G^1])}\big)\\ @AAF^2_{(\w,m,u)}\big([\G^{{1}}]\big) A\\ F^1_{(\w,m,u)}\big(\d_d([\G^1])\big) \end{CD}$$
{given any two elements $[\G^1]$ and $[\G^2] $ of $\P_2(M,*)$, with $\d_d(\G^2)=\d_u(\G^1).$}

{If we are  given elements $[\G^1]$ and $[\G^2] $ of $\P_2(M,*)$, such that $\d_d([\G^2])=\d_u([\G^1]),$ then, for  some representatives $\G_1^1 \in [\G^1]$ and $\G_1^2 \in [\G^2])$ we have $\d_d(\G^2_1)=\d_u(\G^1_1).$ For example, take a rank-1 homotopy $H$ connecting  $\d_u(\G^1)$ with $\d_d(\G^2)$,  and put $\G^1_1=\G^1 \circ_v H$ and $\G^2_1=\G^2$.
Hence $F_{(\w,m,u)}\colon \Pc_2(M,*) \to \C(\Gc)$ is a functor.}
\end{Proof}

{\bf Claim:}
$F_{(\w,m,u)}$ is a monoidal functor.

\begin{Proof}
Let $\G^1$ and $\G^2$ be 2-paths in $\P_2(M,*)$. Let also $\G=\G^1\circ_h \G^2$. Put $\G^1(s,t)=\g_s^1(t)$,  $\G^2(s,t)=\g^2_s(t)$ and $\G(s,t)=\g_s(t)$. Here $s,t \in [0,1]$. We want to prove that $e_\G(s)=\left (g_{\g^1_0}^{-1} \tr {e_{\G^2}(s)}\right) e_{\G^1}(s) $ for any $s \in [0,1]$; see \ref{CCG}. Note that this equation holds for $s=0$. Furthermore we have: 
\begin{align*}
 \frac{d}{d s}& \left (\left (g_{\g^1_0}^{-1} \tr e_{\G^2}(s)\right) e_{\G^1}(s) \right)
\\
&=\Big (  g_{\g^1_0}^{-1} \tr \frac{d}{d s} e_{\G^2}(s)\Big) e_{\G^1}(s) +\left (g_{\g^1_0}^{-1} \tr e_{\G^2}(s)\right)  \frac{d}{d s}  e_{\G^1}(s) 
\\
 &{= \big (g_{\g^1_0}^{-1} \tr \Big (e_{\G^2}(s) \int_{0}^{1} m \left (\t{\frac{\d}{\d t}\g^2_s(t)},\t{\frac{\d}{\d s}\g^2_s(t)} \right)_{{\mathcal{H}}_\w(\g^2_s,t,u) } dt \Big) \Big)e_{\G^1}(s) }
\\
 &\quad \quad \quad \quad+\left (g_{\g^1_0}^{-1} \tr e_{\G^2}(s)\right) e_{\G^1}(s)  \int_{0}^{1} m \left (\t{\frac{\d}{\d t}\g^1_s(t)},\t{\frac{\d}{\d s}\g^1_s(t)} \right)_{{\mathcal{H}}_\w(\g^1_s,t,u) } dt \\&=C_s+{{{D}}}_s.
\end{align*}
{Let us analyse each of the terms above separately:}
\begin{align*}
C_s&=\left (g_{\g^1_0}^{-1} \tr e_{\G^2}(s) \right ) \Big ( g_{\g^1_0}^{-1}\tr \int_{0}^{1} m \left (\t{\frac{\d}{\d t}\g ^2_s(t)},\t{\frac{\d}{\d s}\g^2_s(t)} \right)_{{\mathcal{H}}_\w(\g^2_s,t,u) } dt \Big)e_{\G^1}(s) 
\\
 &=\left (g_{\g^1_0}^{-1} \tr e_{\G^2}(s) \right ) e_{\G^1}(s)\\ &\quad\quad\quad\quad\quad \Big (\big ( \d(e_{\G^1}(s))^{-1} g_{\g^1_0}^{-1} \big) \tr \int_{0}^{1} m \left (\t{\frac{\d}{\d t}\g ^2_s(t)},\t{\frac{\d}{\d s}\g^2_s(t)} \right)_{{\mathcal{H}}_\w(\g^2_s,t,u) } dt \Big) 
\\
&=\left (g_{\g^1_0}^{-1}\tr e_{\G^2}(s) \right)e_{\G^1}(s)\Big ( g_{\g^1_s}^{-1}  \tr \int_{0}^{1} m \left (\t{\frac{\d}{\d t}\g^2_s(t)},\t{\frac{\d}{\d s}\g^2_s(t)} \right)_{{\mathcal{H}}_\w(\g^2_s,t,u) } dt \Big) .
\end{align*}
Given that $m$ is $G$-{equivariant}, given $X,Y \in \X(M)$ we have:
\begin{equation} \label{Refer} 
{g^{-1}\tr m(\t{X},\t{Y})_v=m(\t{X}g,\t{Y}g)_{vg}=m(\t{X},\t{Y})_{vg}.}
\end{equation}
Here $v \in P$ and $g \in G$. Therefore (recall $\g_s=\g^1_s\g^2_s$):
\begin{align*}
C_s&=\left (g_{\g^1_0}^{-1} \tr e_{\G^2}(s) \right)e_{\G^1}(s) \int_{0}^{1} m \left (\t{\frac{\d}{\d t}\g^2_s(t)},\t{\frac{\d}{\d s}\g^2_s(t)} \right)_{{\mathcal{H}}_\w(\g^2_s,t,u)  g_{\g^1_s}  } dt  
\\
 &=\left (g_{\g^1_0}^{-1} \tr e_{\G^2}(s) \right) e_{\G^1}(s)\int_{0}^{1} m \left (\t{\frac{\d}{\d t}\g^2_s(t)},\t{\frac{\d}{\d s}\g^2_s(t)} \right)_{{\mathcal{H}}_\w(\g^2_s,t,ug_{\g^1_s}  )  } dt  
\\
 &=\left (g_{\g^1_0}^{-1} \tr e_{\G^2}(s) \right) e_{\G^1}(s)\int_{1/2}^{1} m \left (\t{\frac{\d}{\d t}\g_s(t)},\t{\frac{\d}{\d s}\g_s(t)} \right)_{{\mathcal{H}}_\w(\g_s,t,u  )  } dt .
\end{align*}
On the other hand
$${{{D}}}_s=\left (g_{\g^1_0}^{-1} \tr e_{\G^2}(s)\right) e_{\G^1}(s)  \int_{0}^{1/2} m \left (\t{\frac{\d}{\d t}\g_s(t)},\t{\frac{\d}{\d s}\g_s(t)} \right)_{{\mathcal{H}}_\w({\g_s,t,u}) } dt. $$
Putting everything together it follows that:
\begin{multline*}
 \frac{d}{d s} \left (\left (g_{\g^1_0}^{-1} \tr {e_{\G^2}(s)}\right) e_{\G^1}(s) \right) \\   =\left (g_{\g^1_0}^{-1} \tr e_{\G^2(s)}\right) e_{\G^1}(s)\int_{0}^{1} m \left (\t{\frac{\d}{\d t}\g_s(t)},\t{\frac{\d}{\d s}\g_s(t)} \right)_{{\mathcal{H}}_\w({\g_s,t,u}) } dt .
\end{multline*}
Therefore $s \in [0,1] \mapsto \left (g_{\g^1_0}^{-1} \tr e_{\G^2(s)}\right) e_{\G^1}(s) \in E$ and $s \in [0,1] \mapsto  e_\G(s) \in E$  are both solutions of the same differential equation  in $E$ and have the same initial condition.

In particular $ \left (g_{\g^1_0}^{-1} \tr e_{\G^2 }\right) e_{\G^1}  = e_\G $. Therefore:
$${
 F\left ( \begin{CD} \d_u([\G^1]) \\ @AA [\G^1]A \\   \d_d([\G^1]) \end{CD} \tn \begin{CD} \d_u([\G^2]) \\ @AA [\G^2]A \\   \d_d([\G ^2]) \end{CD} \right)=\begin{CD} g_{\g^1_1}^{-1} g_{\g^2_1}^{-1} \\ @AA \left (g_{\g^1_0}^{-1} \tr e_{\G^2 }\right) e_{\G^1} A \\g_{\g^1_0}^{-1}g_{\g^2_0}^{-1}  \end{CD}\quad \quad \quad \quad  =\begin{CD} F_1\big(\d_u([\G^1])\big) \\ @AA F_2([\G^1])A \\   F_1\big(\d_d([\G^1])\big) \end{CD} \quad \tn \begin{CD} F_1\big(\d_u([\G^2])\big) \\ @AA F_2\big([\G^2]\big)A \\   F_1\big(\d_d([\G^2])\big )\end{CD} \quad.}
$$

We have therefore proven that 
$F_{(\w,m,u)}\colon \Pc_2(M,*) \to \C(\Gc)$ is a monoidal functor, and therefore a categorical group map.
\end{Proof}

{\bf Claim:}
The functor $F_{(\w,m,u)}\colon \Pc_2(M,*) \to \C(\Gc)$ is smooth both on the sets $\pi_1^1(M,*)$ and $\P_2(M,*)$ of objects and morphisms of $\Pc_2(M,*)$.

\begin{Proof}
This follows directly from the definition of $F_{(\w,m,u)}$, and from the lemma on page 74 of \cite{KN}.
\end{Proof}

The proof of theorem \ref{FFF} is complete.

\subsubsection{{2-Bundles}}\label{2bundles}

Let ${\Gc=(\d\colon E \to G,\tr)}$ be a Lie crossed module. Let  the associated differential crossed module be ${\mathfrak{G}=(\d \colon \le \to  \lg,\tr )}$. Fix a smooth manifold $M$ and a smooth principal $G$-bundle $\pi\colon P \to M$.

A natural weakening of the notion of a categorical connection $(\w,m)$ in $P$ is to consider  $m$ to be only locally defined on $P$,  but so that $\d(m_i)=\W$, for any $i$. Here $m_i \in \A^2(P_i,\le)$, where $P_i=\pi^{-1}(U_i)$ and  $\{U_i\}$ is an open cover of $P$. This can be done by using the general theory of 2-bundles \cite{BS}, or non-abelian gerbes \cite{BrMe}, which correspond to 2-bundles with structure crossed modules of the form {$G=\rm{Aut}(E)$} of example \ref{q}.

To this end, let $\{U_i\}$ be an open cover of $M$. We can suppose that, for any $i_1,\ldots ,i_n$, the $n$-fold intersection $U_{i_1,\ldots ,i_n}\doteq U_{i_1} \cap \ldots \cap U_{i_n}$ is contractible. A $\Gc$-2-bundle can be described by using smooth maps $g_{ij}\colon U_{ij} \to G$ and $h_{ijk}\colon U_{ijk} \to E.$ These {have to satisfy} the cocycle conditions: $\d(h_{ijk}^{-1})g_{ij}g_{jk}=g_{ik}$ and $h_{ijk}h_{ikl}=\left (g_{ij}\tr h_{jkl}\right)h_{ijl}$ for any $i,j,k,l$. Therefore a $\Gc$-2-bundle does not necessarily define a $G$-principal bundle, unless $\d(h_{ijk})=1$ for each $i,j,k$. If this condition is {satisfied}, a $\Gc$-2-bundle will be called a {\it special $\Gc$-2-bundle.}

The main {purpose of} this article was the construction of categorical holonomies, in other words of smooth categorical group maps $\Pc_2(M,*) \to \C(\Gc)$, or, equivalently of smooth crossed module maps {$\Pi_2^{2,1}(M,*)\to \Gc$;} see \ref{thinfund} for this notation.
From the point of view of 2-dimensional holonomy, $\Gc$-2-bundles are natural objects; see \cite{BS}. However,  our notion of categorical holonomy requires the existence of a smooth group morphism $\pi^1_1(M,*) \to G$, and therefore of a principal $G$-bundle with connection; {see \cite{CP}}. This makes it natural to work with special $\Gc$-2-bundles.

A special $\Gc$-2-bundle with connection can be specified by the following data: (Here $A=\ker (\d) \subset E$, and $\la$ is its Lie algebra; note that $A$ is {central in $E$)}.
\begin{enumerate}
 \item A principal $G$-bundle $\pi\colon P \to M$ with a connection $\w$.
\item An open cover $\{U_i\}$ of $M$ such that for any $i_1,\ldots ,i_n$ the $n$-fold intersection $U_{i_1,\ldots ,i_n}$ is contractible. Let $P_{i_1,\ldots ,i_n}=\pi^{-1}(U_{i_1,\ldots ,i_n})$.
\item $G$-{equivariant} $A$-valued functions $e_{ijk}\colon P_{ijk} \to A$. These are required to satisfy:  $e_{ijk}e_{ikl}= e_{jkl} e_{ijl}$.
\item $G$-{equivariant} horizontal 2-forms $m_i\in  \A^2(P_i,\le)$. These are required to satisfy $\d(m_i)=\W$, the curvature 2-form of $\w$.
\item $G$-{equivariant} horizontal 1-forms $\eta_{ij} \in \A^1(P_{ij},\la)$. These are  required to satisfy: $$m_i-m_j=D \eta_{ij} \quad \textrm{ and }\quad  \left(De_{ijk}^{-1}\right) e_{ijk}=\eta_{ij}+\eta_{jk}-\eta_{ik},$$
where $D$ denotes {the} exterior covariant derivative with respect to $\w$, in other words, by using lemma \ref{ext},
$$m_i-m_j=d \eta_{ij} +\w \wedge^\tr \eta_{ij} \quad \textrm{ and }\quad  \left(d e_{ijk}^{-1}+\w \wedge^\tr e_{ijk}^{-1}\right) e_{ijk}=\eta_{ij}+\eta_{jk}-\eta_{ik}.$$
\end{enumerate}
One passes from this setting to the setting in \cite{BS} by choosing local sections $f_i$ of $P_i$, and defining the functions  $g_{ij}\colon U_{ij} \to G$ in the usual way: $f_ig_{ij}=f_j.$ {Next} define $h_{ijk}=f_i^*(e_{ijk})$, $A_i=f_i^*(\w)$, $B_i=f_i^*(m_i)$ and $a_{ij}=f_i^*(\eta_{ij})$. Then this data defines {the local differential forms of Baez and Schreiber \cite{BS}, satisfying equations corresponding to our special $\Gc$-2-bundles with connection, in particular with the {transition functions} $g_{ij}$ satisfying the usual cocycle condition for a principal $G$-bundle. Note that the Baez-Schreiber local differential forms can be seen as a generalisation to arbitrary $\Gc$ of the local description of non-abelian gerbes, with $G=\rm{Aut}(E)$, of Breen and Messing \cite{BrMe}.}

Special $\Gc$-2-bundles with connection define group morphisms {$\pi_1^1(M,*) \to G$.} However, the definition of the categorical holonomy on a {2-track $[\G]\in \P_2(M,*)$} is tricky if one is to deal with the fact that $m$ is not globally defined. This can be done (and will appear in a future work), however, it will only be, {a priori}, defined up {to} multiplication by central elements of $E$, and thus it should not define holonomies $\Pc_2(M,*) \to \C(\Gc)$, unless we restrict to elements of $\pi_2^2(M,*)$, as in \cite{MP}. For instance, it is easy to see that given two categorical $\Gc$-connections $(\w,m)$ and $(\w,m')$ in $P$ then given a {$2$-path} $\G\in \Pc_2(M,*)$ we have
$${e_\G'=\exp  
\Big (
\int_{0}^{1} \eta 
\left (\t{\frac{d}{d t}\g_1(t)} \right)
_{{\mathcal{H}}_\w(\g_1,t,u) }  dt 
 \Big)\exp  
\Big (-
\int_{0}^{1} \eta 
\left (\t{\frac{d}{d t}\g_2(t)} \right)
_{{\mathcal{H}}_\w(\g_2,t,u) }  dt 
 \Big) e_\G}$$
where $m-m'=D \eta=d \eta+\w \wedge^\tr \eta$, moreover {$\eta\in \A^1(P,\la)$} is $G$-{equivariant} and horizontal. We have put $\g_1=\d_u(\G)$ and $\g_2=\d_d(\G)$. Also $e_\G$ and $e_\G'$ denote 2-dimensional holonomy with respect to $(\w,m)$ and $(\w,m')$, respectively.

For this reason we feel that it is likely that the construction in this article, {with a globally defined $m$}, will {describe} all categorical holonomies $\Pc_2(M,*) \to \C(\Gc)$, especially if we consider a generalised unbased setting. We will investigate this issue in a future publication. 
\subsection{Wilson Spheres}

{Let ${\Gc=(\d\colon E \to G,\tr)}$ be a Lie crossed module, and let ${\mathfrak{G}=(\d \colon \le \to  \lg,\tr )}$ be the associated differential crossed module. Fix a smooth manifold $M$ and a smooth principal $G$-bundle $P \to M$, with a connection $\w \in \A^1(P,\lg)$, as well as a $\Gc$-categorical connection $(\w,m)$.  The aim of this subsection is to define the categorical holonomy of a 2-sphere $\Sigma$ embedded in $M$.}

\subsubsection{The dependence of a categorical connection holonomy on the elements of the fibre }

Choose $u,v \in P$. We want to relate $F_{(\w,m,u)}$ with $F_{(\w,m,v)}$. Suppose first that $u,v$ each belong to the fibre $P_x$ of $P$ at a certain point $x\in M$. Let $g\in G$ be the unique element such that $u=vg$. Let $\g$ be a 1-path of $M$ starting and ending at $x$. Then ${\mathcal{H}}_\w(\g,1,u)={\mathcal{H}}_\w(\g,1,vg)={\mathcal{H}}_\w(\g,1,v)g$. Therefore since {$ \mathcal{H}_\w(\g,1,u)F_{(\w,u)}^1(\g)=u$ and ${\mathcal{H}}_\w(\g,1,v)F_{(\w,v)}^1(\g)=v$}
we have the well-known formula:
\begin{equation}\label{A1}F^1_{(\w,v)}(\g)=g F^1_{(\w,u)}(\g) g^{-1}, \end{equation}
for $\g$ a smooth curve starting and ending at $x$ and $u=vg$, where $u,v \in P_x$ and  $g \in G$. 

Let us now see how $F^2_{(m,\w,u)}$  is related to $F^2_{(m,\w,v)}$.
\begin{Lemma}
For any 2-path $\G$ in $M$, with $\d_l(\G)=\d_r(\G)=x$, we have:

 \begin{equation}\label{A2}F^2_{(\w,m,v)}(\G)=g \tr F^2_{(\w,m,u)}(\G), \end{equation}
where $u,v$ belong to the  fibre $P_x$ of $P$, and $g \in G$ is such that  $u=vg$. 
\end{Lemma}
\begin{Proof}
As usual, put $\G(t,s)=\g_s(t); s,t \in [0,1].$
Let $e_\G(s)$ and $f_\G(s)$ be defined by the differential equations:
\begin{equation*}
{\frac{d}{d s} e_\G(s)=e_\G(s)\int_{0}^{1} m \left (\t{\frac{\d}{\d t}\g_s(t)},\t{\frac{\d}{\d s}\g_s(t)} \right)  _{{\mathcal{H}}_\w(\g_s,t,u) } dt,}
\end{equation*}
and
\begin{equation*}
{\frac{d}{d s} f_\G(s)=f_\G(s)\int_{0}^{1} m \left (\t{\frac{\d}{\d t}\g_s(t)},\t{\frac{\d}{\d s}\g_s(t)} \right)  _{{\mathcal{H}}_\w(\g_s,t,v) } dt };
\end{equation*}
with initial conditions  $e_\G(0),f_\G(0)=1_E$; see equation \ref{G}. Then:
\begin{align*}
\frac{d}{d s} g^{-1}\tr f_\G(s)
&=g^{-1} \tr \frac{d}{d s} f_\G(s)\\
&{= g^{-1}\tr  \Big( f_\G(s)  \int_{0}^{1} m \left (\t{\frac{\d}{\d t}\g_s(t)},\t{\frac{\d}{\d s}\g_s(t)} \right)  _{{\mathcal{H}}_\w(\g_s,t,v) } dt \Big)}\\
&{=\left (g^{-1}\tr f_\G(s) \right) \Big(g^{-1}\tr \int_{0}^{1} m \left (\t{\frac{\d}{\d t}\g_s(t)},\t{\frac{\d}{\d s}\g_s(t)} \right) _{{\mathcal{H}}_\w(\g_s,t,v) } dt \Big)}\\
&=\left (g^{-1}\tr f_\G(s) \right)  \int_{0}^{1} m \left (\t{\frac{\d}{\d t}\g_s(t)},\t{\frac{\d}{\d s}\g_s(t)} \right)  _{{\mathcal{H}}_\w(\g_s,t,ug^{-1}) g} dt\\
&=\left (g^{-1}\tr f_\G(s) \right) \int_{0}^{1} m \left (\t{\frac{\d}{\d t}\g_s(t)},\t{\frac{\d}{\d s}\g_s(t)} \right)  _{{\mathcal{H}}_\w(\g_s,t,u) } dt. 
\end{align*}
The fourth equation follows from (\ref{Refer}).
Therefore $g^{-1} \tr f_\G(s)= e_\G(s)$, for each $s \in [0,1]$, thus, in particular,  $F^2_{(\w,m,v)}(\G)=g \tr F^2_{(\w,m,u)}(\G)$.  
\end{Proof}

There is a natural left action of $G$ on {$G\ltimes E$}, given by $g\tr (h,e)=(ghg^{-1},g\tr e)$. Using this, 
we can summarise the results of this subsection, equations (\ref{A1}) and (\ref{A2}), in the following:
\begin{Theorem}\label{I}
If {$v \in P_x$} then
\begin{equation}
{F_{(\w,m,vg)}=g^{-1} \tr F_{(\w,m,v)}}
\end{equation}
for any $g \in G$.
\end{Theorem}

\subsubsection{The dependence of a categorical connection holonomy on the point in the principal bundle}

Let ${\Gc=(\d\colon E \to G,\tr)}$  be a Lie crossed module. Fix a smooth manifold $M$, and a smooth principal $G$-bundle $P \to M$ with a connection $\w \in \A^1(P,\lg)$. Let $(\w,m)$ be a $\Gc$-categorical connection on $P$.
Let $x,y\in M$. Choose $v \in P_y$. Let $\n$ be a 1-path connecting $y$ with $x$. This curve $\n$ defines a categorical group isomorphism $T_\n\colon \Pc_2(M,y)  \to  \Pc_2(M,x)$ in the obvious way, by considering horizontal compositions. More precisely let $\G_\n$ be the 2-path such that $\G_\n(t,s)=\n(t);\forall s,t \in [0,1]$. Let also $\G_{\n^{-1}}(t,s)=\n^{-1}(t);\forall s,t \in [0,1].$
Then {$T_\n([\G])\doteq [\G_{\n^{-1}}] \circ_h [\G]\circ_h [\G_{\n}] \in \P_2(M,x) $} for each $[\G] \in \P_2(M,y),$ and {$T_\n([\g])=[\n^{-1}\g\n]$ for $[\n] \in \pi_1^1(M,y)$}.

\begin{Theorem}
Let {$u={\mathcal{H}}_\w(\n,1,v)$.} We have: 
\begin{equation}
F_{(\w,m,v)}=F_{(\w,m,u)}\circ T_\n.
\end{equation}
\end{Theorem}
\begin{Proof}
Let $[\G] \in \P_2(M,y)$. Then $[\G]=[\G_{\n}] \circ_h  T_\n([\G]) \circ_h [\G_{\n^{-1}}]$. Let $\G^1= \G_{\n} \circ_h \G_{\n^{-1}} \circ_h \G  \circ_h \G_{\n} \circ_h  \G_{\n^{-1}}$, thus $[\G^1]=[\G]$. Put $\g_s(t)=\G^1(t,s); \forall s,t \in [0,1]$.  Let also  $\G'= \G_{\n^{-1}}  \circ_h \G \circ_h \G_{\n}$, thus $[\G']=T_\n[\G].$ Put  $\g_s'(t)=\G'(t,s); \forall s,t \in [0,1]$.
We have
\begin{align*}
\frac{d}{d s} e_{\G^1}(s)&=e_{\G^1}(s)\int_{0}^{1} m \left (\t{\frac{\d}{\d t}\g_s(t)},\t{\frac{\d}{\d s}\g_s(t)} \right)  _{{\mathcal{H}}_\w(\g_s,t,v) } dt\\
&=e_{\G^1}(s)\int_{0}^{1} m \left (\t{\frac{\d}{\d t}\g'_s(t)},\t{\frac{\d}{\d s}\g'_s(t)} \right)  _{{\mathcal{H}}_\w(\g'_s,t,{\mathcal{H}}_\w(\n,1,v)) } dt \\
&{=e_{\G^1}(s)\int_{0}^{1} m \left (\t{\frac{\d}{\d t}\g'_s(t)},\t{\frac{\d}{\d s}\g'_s(t)} \right)  _{{\mathcal{H}}_\w(\g'_s,t,u) } dt }
\end{align*}
{The second  equation follows} since $\G_\n(s,t)$ does not depend on $s$. Comparing with equation \ref{G} {it}  follows that $e_{\G^1}(s)=e_{\G'}(s)$, for each $s \in [0,1]$. Therefore:
\begin{align*}
F_{(\w,m,v)}^2([\G])&=F_{(\w,m,v)}^2([\G^1])=e_{\G^1} =e_{\G'} 
                 =F_{(\w,m,u)}^2([\G'])
                 =F_{(\w,m,u)}^2(T_\n([\G])).
\end{align*}
 This
finishes the proof of the result at the level of the elements of $\P_2(M,y)$. The proof at the level of the  elements of $\pi_1^1(M,y)$ is similar.
 \end{Proof}
\subsubsection{Associating 2-holonomies to embedded spheres}\label{WS}

Let ${\Gc=(\d\colon E \to G,\tr)}$  be a Lie crossed module. Let $M$ be a smooth manifold. Consider a principal $G$-bundle $P$ over $M$ with a $\Gc$-categorical connection $(\w,m)$. Consider an embedded oriented 2-sphere  $\Sigma\subset M$. We want to define the categorical holonomy of $\Sigma${, which we call a {\it Wilson Sphere}, by analogy with Wilson loops, familiar in gauge theory. }

\begin{Definition}
Choose an orientation preserving parametrisation 
$$S^2 = {{\rm D}}^2 / \d {{\rm D}}^2 \ra{\f} \Sigma \subset M$$
 of $\Sigma$. Let $*=\f(\d {{\rm D}}^2) \in \Sigma \subset M$. Choose {$u\in P_*$,} the fibre of $P$ at $*$. The {categorical holonomy} of $\Sigma$ is defined as:
$${\Wc}_{(\w,m)}(\Sigma,\f,u)=F^2_{(\w,m,u)}(\f) \in \ker(\d) \subset E.$$
\end{Definition}
The fact that the {categorical holonomy} ${\Wc}_{(\w,m)}(\Sigma,\f,u)$ lives in $\ker(\d) \subset E$ follows from the fact that {$\d_d(\f)=\d_u(\f)=*$,} and lemma \ref{DDD}.

\begin{Theorem}\label{asd}
The {categorical holonomy} ${\Wc}_{(\w,m)}(\Sigma,\f,u)$ does not depend on $\f$ and $u$, up to acting by an element of $G$. Therefore, we can define a {categorical holonomy} ${\Wc}_{(\w,m)}(\Sigma)$ of embedded 2-spheres $\Sigma$ in $M$ {(the Wilson Sphere)}, taking values in $E/G$. If $\Sigma^*$ is obtained from $\Sigma$ by reversing its orientation then ${\Wc}_{(\w,m)}(\Sigma^*)=\left ({\Wc}_{(\w,m)}(\Sigma)\right)^{-1}$. 
\end{Theorem}

\begin{Example}
{Consider the $U(1)$-bundle $P^n$ over $S^2$ with Chern class $n$, where $n \in \Z$. Define a crossed module $\Gc$ given by the exponential map $\exp \colon i \R \to  U(1)$, considering the trivial action of $U(1)$ on $i\R$. Consider {some} connection 1-form $\w$ on $P_n$, thus $(\w,\W)$ is a categorical $\Gc$-connection on $P^n$. We will then have that $\Wc_{(\w,\W)}(S^2)= n$.}
\end{Example}

\begin{Proof} {\bf (Theorem \ref{asd})}
The independence of ${\Wc}_{(\w,m)}(\Sigma,\f,u)$ on the choice of $u \in P_*$, up to acting by an element of $G$, follows from theorem \ref{I}.

 Choose a different orientation preserving parametrisation  {${{\rm D}}^2 / \d {{\rm D}}^2 \ra{\f'} \Sigma$ of $\Sigma$.}  Since both $\f$ and $\f'$ are orientation preserving, there exists a smooth 1-parameter family of embeddings $x \in [0,1] \mapsto \f_x$ connecting $\f$ and $\f'$. Here $\f_x$ is an orientation preserving parametrisation $\f_x \colon {{\rm D}}^2 / \d {{\rm D}}^2 \to \Sigma$, for each $x \in [0,1]$.

Let $q\colon [0,1] \to M$ be the map such that $q(x)=\f_x(\d ([0,1]^2)), \forall x \in  [0,1]$; therefore $* =q(0)$.  Put  $\G^x(t,s)\doteq \f_x(t,s), $ where $t,s,x \in [0,1]$. Therefore $\d_l(\G^x)=\d_r(\G^x)=q(x),$ for each $x \in  [0,1]$.
 Let also $J(t,s,x)\doteq \G^x(t,s)$, where $t,s,x \in [0,1]$.

The map $J\colon [0,1]^3 \to \Sigma \subset M$ satisfies the conditions of theorem \ref{Main2}. In particular 
$${\Wc}_{(\w,m)}(\Sigma,\f_0,u)=e_{\G^0} =e_{\G^1} ={\Wc}_{(\w,m)}(\Sigma,\f_1,\H_\w(q,1,u)),$$
{since certainly $\Rank ({\EuScript{D}}_z J)\leq 2,\forall z \in [0,1]^3.$}

Finally, let $\f^{-1}(x,y)=\f(1-x,y); \forall x,y \in [0,1]$. It does not follow that $\f^{-1}\circ_h \f$ is rank-2 homotopic to $*$ in {a} strict sense, since a homotopy connecting $\f^{-1}\circ_h \f$ with $*$ will have some singular points, as $\f$ and $\f^{-1}$ are not 2-paths, because they fail the sitting instant conditions.

Let $\g_s(t)=\f(t,s)$ and $\g_s'(t)=\f^{-1}(t,s)$, where $t,s \in[0,1]$. Let us {analyse $e_{\f^{-1}}(s)$ directly.  We have:}
\begin{align*}
\frac{d}{d s} e_{\f^{-1}}(s) 
&=  e_{\f^{-1}}(s)   \int_{0}^{1} m \left (\t{\frac{\d}{\d t}\g_s'(t)},\t{\frac{\d}{\d s}\g_s'(t)} \right)  _{{\mathcal{H}}_\w(\g_s',t,u) } dt\\&= - e_{\f^{-1}}(s)   \int_{0}^{1} m \left (\t{\frac{\d}{\d t}\g_s(t)},\t{\frac{\d}{\d s}\g_s(t)} \right)  _{{\mathcal{H}}_\w(\g_s,t,ug_{\g_s'}) } dt
\end{align*}
Therefore:
\begin{align*}
\frac{d}{d s} e_{\f^{-1}}(s) &= - e_{\f^{-1}}(s) \Big ( g_{\g_s'}^{-1} \tr  \int_{0}^{1} m \left (\t{\frac{\d}{\d t}\g_s(t)},\t{\frac{\d}{\d s}\g_s(t)} \right)  _{{\mathcal{H}}_\w(\g_s,t,u) } dt\Big)\\
&= - e_{\f^{-1}}(s) \Big ( \d \left (e^{-1}_{\f^{-1}}(s) \right ) \tr  \int_{0}^{1} m \left (\t{\frac{\d}{\d t}\g_s(t)},\t{\frac{\d}{\d s}\g_s(t)} \right)  _{{\mathcal{H}}_\w(\g_s,t,u) } dt\Big)\\
&= -\Big( \int_{0}^{1} m \left (\t{\frac{\d}{\d t}\g_s(t)},\t{\frac{\d}{\d s}\g_s(t)} \right)  _{{\mathcal{H}}_\w(\g_s,t,u) } dt\Big)  e_{\f^{-1}}(s).
\end{align*}
The result follows from the fact that if $A(t)$ is a smooth function $[0,1] \to E$ and if $w(s)$ is a solution of $\frac{d}{d s} w(s)=w(s)A(s)$ then the smooth function  $v(s)\doteq w(s)^{-1}$ is a solution of $\frac{d}{d s} v(s)=-A(s)v(s)$.  In particular $ e_{\f^{-1}}(1)=  e_{\f}^{-1}(1)$.
\end{Proof}

Note that the previous proof uses, in an essential way, the fact that the mapping class group of $S^2$ is $\{\pm 1\}$. Therefore, it is likely that the previous theorem does not extend to an embedded surface of genus greater than zero. However, it is probable that this can be fixed if any embedded surface $\Sigma \subset M$ is assigned  an equivalence class of embedding $S \to \Sigma$, where two embeddings $\f,\f'\colon S \to \Sigma$ are equivalent if $\f^{-1} \circ \f'$ is isotopic to the identity diffeomorphism. {We will analyse this issue in a subsequent publication. }

\section{Appendix: Technical Lemmas}
The following lemmas provide the foundation for the construction of the 2-category $\Sc_2(M)$, defined for any smooth manifold $M$. Specifically, lemma \ref{MAIN} is used to prove that the vertical composition of 2-tracks is well defined; see \ref{HVC}.

\begin{Lemma}
Let $M$ be a smooth manifold. Let $f \colon S^2 \to M$ be a smooth map. Consider a handle decomposition $H$ of $M$. Let $H^i$ be the $i$-skeleton of $H$; in other words $H^i$ is  the handlebody made from the handles of $H$ of index less than or equal to $i$, where $i=0,1, \ldots ,\dim(M)$. There exists an {ambient isotopy} of $M$, sending $H$ to a handle decomposition $I$ of $M$, such that $f(S^2) \subset I^2$, the 2-skeleton of the handle decomposition $I$.

\end{Lemma}

\begin{Proof}
Let $m=\dim(M)$.

Let $J$ be the handle decomposition of $M$ which is dual to $H$. If $m=2$ there is nothing to prove. Otherwise $f(S^2)$ is a compact set with  zero measure in $M$, by Sard's theorem; see \cite{H}. Therefore, there exists an {ambient isotopy} of $M$, sending $J$ to a handle decomposition  $J_0$ of $M$, such that the 0-handles of $J_0$ are contained in $M \setminus f(S^2)$, a non-empty open set of $M$. If $m=3$, taking the dual handle decomposition of $J_0$, yields  a handle decomposition $I$ of $M$ with $f(S^2) \subset I^2$.

Suppose  that $m>3$. Let $\{h_i^1\}$ be the set of handles of $J_0$ of index $1$. Then each $1$-handle $h_i^1$ is of the form $h_i^1={{\rm D}}^1_i \times {{\rm D}}^{m-1}_i$, and attaches to the $0$-skeleton $J^0_0$ of $J_0$ along $S^0_i \times {{\rm D}}^{m-1}_i$. For each $i$, let $p_i^1 \colon h_i^1 \to {{\rm D}}^{m-1}_i$ be the obvious projection. {Here ${\rm D}^m$ is the $m$-disk.} Then $p_i^1$ {can be extended} to a smooth map $q_i^1 \colon V_i^1 \to \R^{m-1}$, where $V_i^1$ is an open neighbourhood of $h_i^1$ in $M$.  For each $i$, let $A_i^1=f^{-1}(V_i^1)\subset S^2$. It is an open set in $S^2$. Since $m-1 >2$, the set $(q_i^1 \circ f)(A_i^1)\subset \R^{m-1}$ has measure zero, for  $q_i^1$ is smooth. Therefore $p_i^1(f(S^2))\subset (q_i^1 \circ f)(A_i^1)$ has measure zero in ${{\rm D}}^{m-1}_i$.

 Choose an  $(m-1)$-ball $c^{m-1}_i \subset {{\rm D}}^{m-1}_i$, not intersecting $p_i^1(f(S^2))$; note that $p_i^1(f(S^2))$ is compact. There exists an {ambient isotopy} of $M$, which restricts to an {ambient isotopy} of $J^0_0$ (the $0$-skeleton of $J_0$), whilst sending each $h_i^1={{\rm D}}^1_i \times {{\rm D}}^{m-1}_i$ to ${{\rm D}}^1_i\times c^{m-1}_i$. Let $J_1$ be the handle decomposition of $M$  to which $J_0$ is sent by the {ambient isotopy}. By construction $J^1_1 \cap f(S^2) = \emptyset.$ If $\dim(M)=4$, taking the dual handle decomposition to $J_1$, yields the handle decomposition $I$ of $M$, and by construction $ f(S^2) \subset I^2.$

Suppose $m>4$. The argument is the same. As above, let $\{h_i^2\}$ be the set of  2-handles of $J_1$. Then, for any $i$, $h^2_i={{\rm D}}^2_i \times {{\rm D}}^{m-2}_i$,  and $h^2_i$ attaches to $J_1^1$ along $S^1_i \times {{\rm D}}^{m-2}_i$.  For each $i$, let $p_i^2 \colon h_i^2 \to {{\rm D}}^{m-2}_i$ be the obvious projection. As before, since $m-2 >2$, the set $p_i^2(f(S^2))$, which is compact,  has zero  measure in ${{\rm D}}^{m-2}_i$. Therefore, there exists an  $(m-2)$-ball $c^{m-2}_i \subset {{\rm D}}^{m-2}_i$ such that  $p_i^2(f(S^2)) \cap c^{m-2}_i  = \emptyset$. There exists an {ambient isotopy} of $M$, restricting to an {ambient isotopy} of  $J^1_1$, and sending $h_i^2={{\rm D}}^2_i \times {{\rm D}}^{m-2}_i$ to  ${{\rm D}}^2_i\times c^{m-2}_i$, for each $i$. Let the handle decomposition $J_2$ be the image of $J_1$ under the {ambient isotopy}. As before, by construction,  $J^2_2 \cap f(S^2) = \emptyset$. If $m=5$, letting $I$ be the handle decomposition of $M$ dual to $J_2$, then we have that $f(S^2) \subset J^2$.

{An obvious inductive argument will finish the proof for any $m \in \N$.} \end{Proof}

\begin{Lemma}\label{MAIN}

Let $f \colon \d ({{\rm D}}^3) \to M$ be a smooth map such that $\Rank ({{\EuScript{D}}}_v f)\leq 1, \forall v \in \d ({{\rm D}}^3)$. Here ${{\rm D}}^3=[0,1]^3$. Suppose that $f$ is constant in a neighbourhood of each vertex of $\d({{\rm D}}^3) $. In addition, suppose also that in a neighbourhood $I \times [-\epsilon,\epsilon]$ of each edge $I$ of $\d ({{\rm D}}^3)$, {$f(x,t)=\f(x)$, where $(x,t) \in I \times [-\epsilon,\epsilon]$ and $\f\colon I \to M$ is smooth.}
 Then $f$ can be extended to a smooth map $F \colon {{\rm D}}^3 \to M$ such that  $\Rank ({{\EuScript{D}}}_w F)\leq 2, \forall w \in {{\rm D}}^3$. Moreover we can choose  $F$  so that it has a product structure close to the boundary of ${{\rm D}}^3$.
\end{Lemma}

\begin{Proof}
The conditions of the statement imply that $f$ is of the form $f' \circ r$, where $f'\colon S^2 \to M$ is a smooth map and {$r \colon \d ({{\rm D}}^3) \to S^2$} is some smooth homeomorphism. 

By using  the previous lemma,  we can see that there exists a handle decomposition of $M$ such that $f(  \d ({{\rm D}}^3)  )$ is contained in the 2-skeleton $M^2$ of $M$ (the handlebody made from the 0-,1- and 2-handles of $M$). Let us prove that we can suppose further that $f( \d ({{\rm D}}^3) )$ is contained in the 1-skeleton $M^1$ of $M$. This will be completely analogous  (though dual) to the proof of the previous lemma.

Let $n=\dim(M)$. Suppose that $n \ge 2$. Let $\{h_i\}$ be the set of 2-handles of $M$. Therefore, for each $i$,  $h_i={{\rm D}}^2_i \times {{\rm D}}^{n-2}_i$,  and $h_i$ attaches to $M^1$ along $S^1_i \times {{\rm D}}^{n-2}_i=(\d {{\rm D}}^2_i) \times {{\rm D}}^{n-2}_i$. It suffices to prove that for each $i$ there exists  an $x_i \in {{\rm D}}^2_i$ (thus an open 2-ball) such that $f(\d {{\rm D}}^3)$ does not intersect $\{x_i\} \times {{\rm D}}^{n-2}_i$.

Let $p_i \colon h_i={{\rm D}}^2_i \times {{\rm D}}^{n-2}_i \to {{\rm D}}^2 \times \{0\}$ be the obvious projection. {Given that} $\Rank({{\EuScript{D}}}_v f) \leq 1, \forall v \in  {\d {{\rm D}}^3}$, it follows that $Z_i=p_i(h_i \cap f(\d {{\rm D}}^3))$ has zero measure in ${{\rm D}}^2$, thus that we can choose an $x_i$ in  the interior of ${{\rm D}}^2$ not intersecting $Z_i$.  {Therefore, some open ball in ${{\rm D}}^2$} containing $x_i$ does not intersect $Z_i$. This proves that there exists a handle decomposition of $M$ such that  $f( \d {{\rm D}}^3  )$ is contained in the interior of the 1-skeleton $M^1$ of $M$,  {analogously to the proof of the previous lemma.}

Consider the universal cover  $N \ra{p} M^1$ of $M^1$. It is a contractible smooth manifold.  We can lift $f$ as $f=p \circ {g'}$, where ${g'}\colon \d {{\rm D}}^3  \to N$ is a smooth map, such that $\Rank ({{\EuScript{D}}}_v {g'})\leq 1, \forall v \in  {\d{{{\rm D}}^3}}${, since $p$ is a local diffeomorphism}. Let $c\colon N \times [0,1] \to N$ be a smooth contraction of $N$ to a point of it. We can suppose that there exists an $\e >0$ such that $c$ does not depend on $t$ for $t \in [0,\e] \cup [1-\e,1]$. The smooth homotopy $g(v,t)=p ( c({g'}(v),t))$ where $t \in [0,1]$ and $v \in   \d ({{\rm D}}^3) $ will permit us to extend $f \colon  \d ({{\rm D}}^3)  \to M$ to a smooth map $F\colon {{\rm D}}^3 \to M$, with a product structure close to the boundary of ${{\rm D}}^3$. By construction we also have that $\Rank ({{\EuScript{D}}}_w F)\leq 2, \forall w \in {{\rm D}}^3$.
\end{Proof}

\end{document}